\documentclass[reqno, 11pt, letterpaper]{amsart}

\usepackage{amsmath}
\usepackage{amsfonts}
\usepackage{amssymb}
\usepackage{graphicx}
\usepackage{amsthm,graphicx,color,yfonts, mathtools, enumerate}
\usepackage{hyperref}
\usepackage{mathabx}
\usepackage{array}

\usepackage{physics}

\newtheorem{theorem}{Theorem}

\newtheorem{proposition}[theorem]{Proposition}
\newtheorem{lemma}[theorem]{Lemma}

\newtheorem{assumption}[theorem]{Assumption}

\theoremstyle{remark}
\newtheorem{remark}[theorem]{Remark}

\newcommand{\R}{\mathbb{R}}

\usepackage{setspace}
\usepackage{wrapfig}
\usepackage{tikz}
\usetikzlibrary{arrows,calc,decorations.pathreplacing}
\definecolor{light-gray1}{gray}{0.90}
\definecolor{light-gray2}{gray}{0.80}
\definecolor{light-gray3}{gray}{0.60}
\usepackage{bbm}

\makeatletter
\newcommand*{\rom}[1]{\expandafter\@slowromancap\romannumeral #1@}
\makeatother

\numberwithin{equation}{section}

\numberwithin{theorem}{section}

\numberwithin{table}{section}

\numberwithin{figure}{section}

\ifx\pdfoutput\undefined
  \DeclareGraphicsExtensions{.pstex, .eps}
\else
  \ifx\pdfoutput\relax
    \DeclareGraphicsExtensions{.pstex, .eps}
  \else
    \ifnum\pdfoutput>0
      \DeclareGraphicsExtensions{.pdf}
    \else
      \DeclareGraphicsExtensions{.pstex, .eps}
    \fi
  \fi
\fi

\title[Semiclassical equivalence of white dwarf models]{Semiclassical equivalence of two white dwarf models as ground states of the relativistic Hartree-Fock and Vlasov-Poisson energies}

\date{\today}
\author[Y. Hong]{Younghun Hong}
\address{Department of Mathematics, Chung-Ang University, Seoul 06974, Korea}
\email{yhhong@cau.ac.kr}

\author[S. Jin]{Sangdon Jin}
\address{Department of Mathematics Education, Chungbuk National University,  Cheongju 28644,  Korea}
\email{sangdonjin@cbnu.ac.kr}

\author[J. Seok]{Jinmyoung Seok}
\address{Department of Mathematics Education, Seoul National University, Seoul 08826, Korea}
\email{jmseok@snu.ac.kr}

\begin{document}
\maketitle

\begin{abstract}

We are concerned with the semi-classical limit for ground states of the relativistic Hartree-Fock energies (HF) under a mass constraint, which are considered as the quantum mean-field model of white dwarfs \cite{LeLe}. 
In Jang and Seok \cite{JS}, fermionic ground states of the relativistic Vlasov-Poisson energy (VP) are constructed as a classical mean-field model of white dwarfs, and are shown to be equivalent to the classical Chandrasekhar model. In this paper, we prove that as the reduced Planck constant $\hbar$ goes to the zero, the $\hbar$-parameter family of the ground energies and states of (HF) converges to the fermionic ground energy and state of (VP) with the same mass constraint. 
\end{abstract}

\section{Introduction}
\subsection{Backgrounds and setup}
White dwarfs are compact and dense celestial objects supporting themselves by the electron degeneracy pressure against the gravitational collapse. 
They are considered as final evolutionary states of relatively light stars that have exhausted the nuclear fuel.
In 1931, in his seminal paper \cite{C}, Chandrasekhar discovered the existence of the maximal mass $M_c$, the so-called Chandrasekhar limit mass,  of non-rotating stable white dwarfs so that stars with larger mass collapse to other states such as neutron stars or black holes.

To be precise, the densify function $\rho$ of a non-rotating radially symmetric white dwarf can be described by a solution to the equation for a gravitational hydrodynamic equilibrium 
\begin{equation}\label{CE}
\frac{1}{r^2}\frac{d}{dr}\left(\frac{r^2}{\rho}\frac{dP}{dr}\right) = -4\pi \rho,
\end{equation}
where $r$ denotes the radial variable in $\R^3$, $P$ is the degeneracy pressure of the electron gas and the gravitational constant is normalized.  
In \cite{C}, the local equation of state 
\begin{equation}\label{EOS}
P(\rho) =  \frac{4\pi}{3}\int_0^{(  3\rho/4\pi)^{1/3}} \frac{u^4}{\sqrt{1+u^2}}\,du,
\end{equation}
is derived for white dwarfs. The equation \eqref{CE} equipped with \eqref{EOS} is referred to the Chandrasekhar equation.

From a variational formulation, the Chandrasekhar equation is derived as the Euler-Lagrange equation for an energy minimizer subject to the mass constraint,
\begin{equation}\label{Ch}
\mathcal{E}_\textup{C;min}(M) \coloneqq \inf\left\{\mathcal{E}_\text{C}(\rho) ~\Big|~ \rho \geq 0, \, \rho \in L^{4/3}(\R^3),\, \int_{\mathbb{R}^3}\rho(x) dx = M \right\},
\end{equation}
where  $M>0$ and the Chandrasekhar energy is given by 
\[
\mathcal{E}_\textup{C}(\rho) \coloneqq \underbrace{\int_{\mathbb{R}^3} A(\rho)\,dx}_{\textup{kinetic energy}} \underbrace{- \frac{1}2\iint_{\mathbb{R}^3\times\mathbb{R}^3}\frac{\rho(x)\rho(x')}{|x-x'|}\,dxdx'}_{\textup{potential energy}}
\]
with
$$A(\rho) = \int_0^\rho\big(\sqrt{1+(\tfrac{3}{4\pi}u)^{2/3}}-1\big) du.$$
For this hydrodynamic model, the existence of the critical mass $M_c>0$ is established in Lieb and Yau \cite{LY}; $-\infty < \mathcal{E}_\textup{C;min}(M)< 0$ if $0 < M < M_c$, while $\mathcal{E}_\textup{C;min}(M)=-\infty$ if $M > M_c$. It is also shown that in the former case, $\mathcal{E}_\textup{C;min}(M)$ has a unique radially symmetric compactly supported minimizer $\rho_\text{C}$, and that $\rho_\text{C}$ solves the Chandrasekhar equation.

\subsubsection{Relativistic Hartree-Fock formulation}

For the semi-classical description of relativistic white dwarfs, we follow the mathematical formulation in \cite{FL, HS, HLLS}   with the small parameter $\hbar>0$, which represents the reduced Planck constant. In the Heisenberg picture, a relativistic fermionic gas is described by a compact self-adjoint operator $\gamma$ acting on $L^2(\mathbb{R}^3)$. By the Pauli exclusion principle, we assume that $0\leq \gamma\leq 1$ as a quadratic form. For such an  operator $\gamma$, we denote its integral kernel by $\gamma(x,x')$, i.e., 
$$(\gamma\phi)(x)=\int_{\mathbb{R}^3}\gamma(x,x')\phi(x')dx',$$
and the semi-classical density function by 
$$\rho_\gamma^\hbar \coloneqq (2\pi\hbar)^3\rho_\gamma\textup{ with }\rho_\gamma(x)=\gamma(x,x).$$
We define the quantum self-generated potential by
$$\Phi_{\gamma}^\hbar:=-\frac{1}{|x|}*\rho^\hbar_{\gamma}$$
and the exchange term $X_\gamma^\hbar$ as the operator with kernel
$$X_\gamma^\hbar(x,x'):=\frac{(2\pi\hbar)^3\gamma(x,x')}{|x-x'|}.$$
Using the above definitions, we introduce the semi-classically scaled relativistic Hartree-Fock energy functional
\begin{equation}\label{HF energy}
\begin{aligned}\mathcal{E}_{\textup{HF}}^\hbar(\gamma)& \coloneqq \underbrace{\textup{Tr}^\hbar\left((\sqrt{1-\hbar^2\Delta}-1)\gamma\right)}_{\textup{kinetic energy}}+\underbrace{\frac{1}{2}\textup{Tr}^\hbar\big(\Phi_{\gamma}^\hbar\gamma\big)}_{\textup{potential energy}}+\underbrace{\frac{1}{2}\textup{Tr}^\hbar\big(X_\gamma^\hbar\gamma\big)}_{\textup{exchange energy}}\\
\end{aligned}
\end{equation}
and the mass
$$\mathcal{M}^\hbar(\gamma) \coloneqq \textup{Tr}^\hbar(\gamma),$$
where
$$\textup{Tr}^\hbar(\cdot) \coloneqq (2\pi\hbar)^3\textup{Tr}(\cdot).$$
Note that in the kinetic energy, with an abuse of notation, we write $\textup{Tr}^\hbar(\sqrt{1-\hbar^2\Delta}\gamma)=\textup{Tr}^\hbar((1-\hbar^2\Delta)^{\frac{1}{4}}\gamma(1-\hbar^2\Delta)^{\frac{1}{4}})$, where $(1-\hbar^2\Delta)^{\frac{1}{4}}$ is the Fourier multiplier operator with symbol $(1+\hbar^2|\xi|^2)^{\frac{1}{4}}$, because it holds by cyclicity of the trace when $\gamma$ is regular enough.

We define the quantum mechanical admissible class $\mathcal{A}_{\textup{qm}}$ as the collection of compact self-adjoint operators $\gamma: L^2(\mathbb{R}^3)\to L^2(\mathbb{R}^3)$ such that $0\leq\gamma\leq1$ and $\textup{Tr}\big((1-\Delta)^{\frac{1}{4}}\gamma(1-\Delta)^{\frac{1}{4}}\big)<\infty$. Then, the fermionic gas in a white dwarf with mass $M$ is formulated a minimizer of the variational problem 
\begin{equation}\label{var-HF}
\mathcal{E}_{\textup{HF}; \min}^\hbar(M)=\inf_{\gamma\in\mathcal{A}_{\textup{qm}}(M)}\mathcal{E}_{\textup{HF}}^\hbar(\gamma),
\end{equation}
where  $M>0$ and  $\mathcal{A}_{\textup{qm}}(M) \coloneqq \{ \gamma \in \mathcal{A}_{\textup{qm}} ~|~ \mathcal{M}^\hbar(\gamma)=M\}$. In \cite{LeLe}, it is shown that $\mathcal{E}_{\textup{HF}; \min}^\hbar(M)$ has a minimizer, provided that $M$ is strictly less than the \textit{quantum Chandrasekhar limit mass}
$$\mathbf{M}_{\textup{qm}} \coloneqq (2K_{\textup{qm}})^{\frac{3}{2}},$$
where $K_{\textup{qm}}$ is the sharp constant for the Lieb-Thirring type inequality, i.e.,  
\begin{equation}\label{quantum sharp constant}
K_{\textup{qm}} \coloneqq \inf_{\gamma\in \mathcal{A}_{\textup{qm}}}\frac{\|\gamma\|^{\frac{1}{3}}\big(\textup{Tr}^\hbar(\gamma)\big)^{\frac{2}{3}}\textup{Tr}^\hbar\big(|\hbar\nabla|^{\frac{1}{2}}\gamma |\hbar\nabla|^{\frac{1}{2}}\big)}{\frac{1}{4\pi}\|\nabla \Phi_{\gamma}^\hbar\|_{L^2(\mathbb{R}^3)}^2}.
\end{equation}
For $\mu<0$, we denote the characteristic function on $(-\infty, \mu)$ by
$$\mathbbm{1}_{(E<\mu)}:=\left\{\begin{aligned}
&1&&\textup{if }E<\mu,\\
&0&&\textup{if }E\geq \mu.
\end{aligned}\right.$$

\begin{theorem}[$\textup{\cite[Lenzmann-Lewin]{LeLe}}$]\label{thm: quantum existence}
Let $\hbar>0$. For $0<M<\mathbf{M}_{\textup{qm}}$, the following hold.
\begin{enumerate}[$(i)$]
\item The variational problem $\mathcal{E}_{\textup{HF}; \min}^\hbar(M)$ has a minimizer $\mathcal{Q}_\hbar$.
\item The operator $\sqrt{1-\hbar^2\Delta}-1+\Phi_{\mathcal{Q}_\hbar}^\hbar+X_{\mathcal{Q}_\hbar}^\hbar$ has a negative eigenvalue $\mu_\hbar$ such that 
\begin{equation}\label{eq: quantum self-consistent eq}
\mathcal{Q}_\hbar= \mathbbm{1}_{(\sqrt{1-\hbar^2\Delta}-1+\Phi_{\mathcal{Q}_\hbar}^\hbar+X_{\mathcal{Q}_\hbar}^\hbar<\mu_\hbar)}+\mathcal{R}_\hbar,
\end{equation}
where $\mathcal{R}_\hbar$ is a finite rank operator on the eigenspace corresponding to $\mu_\hbar$.
\end{enumerate}
\end{theorem}

\begin{remark}
\begin{enumerate}
\item In Lenzmann-Lewin \cite{LeLe}, the existence and properties of minimizers are proved for the more delicate Hartree-Fock-Bogoliubov energy, but their approach can be adapted to the Hartree-Fock case.
\item By homogeneity, the right hand side of \eqref{quantum sharp constant} is invariant under the choice of $\hbar>0$, and so is the critical mass $\mathbf{M}_{\textup{qm}}$.
\item By functional calculus, $\mathbbm{1}_{(\sqrt{1-\hbar^2\Delta}-1+\Phi_{\mathcal{Q}_\hbar}^\hbar+X_{\mathcal{Q}_\hbar}^\hbar<\mu_\hbar)}$ is the spectral projection onto the eigenspace corresponding to the eigenvalues of $\gamma$ strictly less than $\mu_\hbar$.
\end{enumerate}
\end{remark}

\subsubsection{Relativistic Vlasov-Poisson formulation}

Next, we review an analogous kinetic description for white dwarfs from \cite{JS}. In kinetic theory, relativistic fermionic gases are described as distributions on the phase space. For a distribution $f=f(q,p):\mathbb{R}^3\times\mathbb{R}^3\to[0,1]$, we define the mass by 
$$\mathcal{M}(f) \coloneqq \iint_{\mathbb{R}^3\times\mathbb{R}^3}f(q,p)dqdp,$$
and the relativistic Vlasov-Poisson energy functional (see \cite{GS, HR, JS, LMR}) by 
$$\mathcal{E}_{\textup{VP}}(f)\coloneqq \underbrace{\iint_{\mathbb{R}^3\times\mathbb{R}^3}\big(\sqrt{1+|p|^2}-1\big)f(q,p)dqdp}_{\textup{kinetic energy}}+\underbrace{\frac{1}{2}\iint_{\mathbb{R}^3\times\mathbb{R}^3}\Phi_f(q) f(q,p)dqdp}_{\textup{potential energy}},$$
where 
$$\Phi_{f}:=-\frac{1}{|x|}*\rho_{f}$$
is the classical self-generated potential given by the kinetic density function
$$\rho_f=\int_{\mathbb{R}^3}f(\cdot,p)dp.$$
We define the \textit{classical mechanical} admissible class by
$$\mathcal{A}_{\textup{cm}} \coloneqq \Big\{f\in L^1(\R^3\times\mathbb{R}^3):\ 0\leq f(q,p)\leq 1\textup{ a.e.,}\textup{ and }|p|f\in L^1(\mathbb{R}^3\times\mathbb{R}^3)\Big\},$$
where the condition $0\leq f\leq 1$ is from the Pauli exclusion principle for fermions. Then, the particle distribution in a white dwarf with mass $M$ is modeled as a minimizer for the classical variational problem 
\begin{equation}\label{var-VP}
\mathcal{E}_{\textup{VP};\min} (M) \coloneqq \inf_{f\in \mathcal{A}_{\textup{cm}}(M)}\mathcal{E}_\textup{VP}(f),
\end{equation}
where $M>0$ and $\mathcal{A}_\textup{cm}(M) \coloneqq \{ f \in \mathcal{A}_\textup{cm} ~|~ \mathcal{M}(f) = M\}$. Recently, in Jang-Seok \cite{JS}, it is constructed when its mass is strictly less than the \textit{classical Chandrasekhar limit mass}. Indeed, such a critical mass is chosen analogously. Precisely, it is given by 
$$\mathbf{M}_{\textup{cm}} \coloneqq (2K_{\textup{cm}})^\frac32,$$
where $K_{\textup{cm}}$ is the sharp constant for the corresponding kinetic interpolation inequality 
\begin{equation}\label{classical sharp constant}
K_{\textup{cm}}:=\inf_{f\in \mathcal{A}_\textup{cm}}\frac{\|f\|_{L^\infty(\R^3\times\mathbb{R}^3)}^{\frac{1}{3}}\|f\|_{L^1(\R^3\times\mathbb{R}^3)}^{\frac{2}{3}}\||p|f\|_{L^1(\R^3\times\mathbb{R}^3)}
}{\frac{1}{4\pi}\|\nabla\Phi_{f}\|_{L^2(\mathbb{R}^3)}^2}.
\end{equation}

\begin{theorem}[$\textup{\cite[Jang-Seok]{JS}}$ ]\label{thm: classical existence}
If $0<M<\mathbf{M}_{\textup{cm}}$, the following hold.
\begin{enumerate}[$(i)$]
\item The variational problem $\mathcal{E}_{\textup{VP}; \min}(M)$ has a minimizer $\mathcal{Q}$, which is  unique up to  a translation.
Moreover, for any minimizing sequence $\{f_n\}$ for $\mathcal{E}_{\textup{VP}; \min}(M)$, one has  $\|\nabla(\Phi_{f_n} -\Phi_{\mathcal{Q}})\|_{L^2(\R^3)} \to 0$ as $n\to \infty$, up to a translation. 
\item For some $\mu<0$, $\mathcal{Q}$ obeys the self-consistent equation
\begin{equation}\label{eq: kinetic self-consistent eq}
\mathcal{Q}= \mathbbm{1}_{(\sqrt{1+|p|^2}-1+ \Phi_\mathcal{Q}(q)<\mu)}.
\end{equation}
\item The density function $\rho_\mathcal{Q}$ attains the minimum Chandrasekhar energy, i.e., $\mathcal{E}_\textup{C}(\rho_\mathcal{Q})=\mathcal{E}_\textup{C;min}(M)$ \eqref{Ch}, and consequently it satisfies the Chandrasekhar equations \eqref{CE} and \eqref{EOS}.
\end{enumerate}
\end{theorem}
 \begin{remark}
In Theorem \ref{thm: classical existence} $(i)$, by the uniqueness, the sequential convergence is obtained (see \cite[Remark 2.7]{JS} and \cite{ LS}).
\end{remark}

\subsection{Statement of the main result}
In this paper, we are interested in the connection between the quantum and the classical variational problems \eqref{var-HF} and \eqref{var-VP} via the semi-classical limit. For the statement, we remark that the quantum critical mass is less than or equal to the classical one; $\mathbf{M}_{\textup{qm}} \leq \mathbf{M}_{\textup{cm}}$ (see Lemma \ref{critical mass comparison}). Our main theorem asserts that if $M$ is strictly less than the quantum critical mass, then the ground energies and minimizers correspond each other as stated below.

\begin{theorem}\label{main-theorem}
Let $M \in (0, \mathbf{M}_{\textup{qm}})$ be given. 
Let $\mathcal{Q}_\hbar$ and $\mathcal{Q}$  be  a family of minimizers and a minimizer of the variational problems  \eqref{var-HF} and \eqref{var-VP} and respectively. 
Then, as $\hbar\to0$, after taking suitable spatial translations, one has
\[
\left\{\begin{aligned}
&\mathcal{E}_{\textup{HF}; \min}^\hbar(M) \to \mathcal{E}_{\textup{VP}; \min}(M) \\
&\rho_{\mathcal{Q}_{\hbar}} \rightharpoonup \rho_{\mathcal{Q}}\text{ weakly in } L^q(\R^3) \quad \forall q > 1 \\
&\| \Phi_{\mathcal{Q}_{\hbar}}  -\Phi_{\mathcal{Q}}\|_{\dot{H}^1(\R^3) \cap C^1(\R^3)} \to 0.
\end{aligned}\right.
\]
\end{theorem}

\begin{remark} 
\begin{enumerate}[$(i)$]

\item In Lieb and Yau \cite{LY, LY1}, Chandrasekhar's hydrodynamic white dwarf model is derived directly from the mean-field limit $(N=\hbar^{-\frac13}\to\infty)$ of the ground energy of the fermionic $N$-body relativistic Schr\"odinger Hamiltonian 
$$H_N = \sum_{i=1}^N\sqrt{-\Delta_{x_i}+1}-1 -\frac{1}{N^{2/3}}\sum_{1\leq i< j\leq N}\frac{1}{|x_i-x_j|}.$$
Our main result, together with Theorem \ref{thm: classical existence} $(iii)$, justifies the consistency among three well-known mean-field models; the Hartree-Fock, the Vlasov-Poisson and the Euler-Poisson energies.

\item Contrary to the mean-field limit \cite{LY, LY1}, in our setting, one can take advantages from that minimizers are solutions to the self-consistent equations \eqref{eq: quantum self-consistent eq} and \eqref{eq: kinetic self-consistent eq}. It allows us to obtain a wider range of weak convergence for density functions and the strong convergence for potential functions.
 
\item In Choi, Hong and Seok \cite{CHS}, a similar analysis is presented for free energy minimizers between the non-relativistic Hartree and the Vlasov-Poisson models, where the free energies are defined as the sum of the energy and a Casimir functional (or a generalized entropy). In this work, the Casimir functional plays a role of preventing the ground energy from diverging to $-\infty$. On the other hand, in the present paper, Pauli's exclusion principle ($0\leq\gamma\leq1$ and $0\leq f\leq 1$) has the same role.

\end{enumerate}
\end{remark}

The proof of the main result is based on the approach of Choi-Hong-Seok \cite{CHS}. 
In our proof, a key step is to show the convergence from the quantum to the classical minimum energy, i.e., 
$\mathcal{E}_{\textup{HF}; \min}^\hbar(M) \to \mathcal{E}_{\textup{VP}; \min}(M)$ as $\hbar\to 0$; then the convergence of the density functions of minimizers  $\mathcal{Q}^\hbar$ and  $\mathcal{Q}$ follows from uniform regularity (see Lemma \ref{regularity lemma})  and  the convergence of ground energy levels between $\mathcal{E}_{\textup{HF}; \min}^\hbar(M)$ and $\mathcal{E}_{\textup{VP}; \min}(M)$.
We notice that this key step cannot be directly achieved because the quantum and the classical minimizers are different mathematical objects 
-- one is an operator and the other is a distribution function on the phase space. So, it is not possible to insert the classical minimizer $\mathcal{Q}$ to the Hartree-Fock energy $\mathcal{E}_{\textup{HF}}^\hbar(\gamma)$ or the quantum minimizer $\mathcal{Q}^\hbar$ to the Vlasov-Poisson energy $\mathcal{E}_{\textup{VP}}$ for energy comparison.

They nevertheless have common structures, that is,   both are characteristic functions of local Hamiltonians (see \eqref{eq: quantum self-consistent eq} and \eqref{eq: kinetic self-consistent eq}). 
Thus, we introduce an auxiliary quantum state $\gamma_\hbar=\mathbbm{1}_{(\sqrt{1-\hbar^2\Delta}-1+\Phi_{\mathcal{Q}}< \tilde{\mu}_\hbar)}+\tilde{\mathcal{R}}_\hbar$ mimicking the form of the classical energy minimizer $\mathcal{Q}=\mathbbm{1}_{(\sqrt{1+|p|^2}-1+\Phi_{\mathcal{Q}}<\mu)}$, where $\Phi_{\mathcal{Q}}$ is the potential function of the classical minimizer $\mathcal{Q}$ (see Lemma \ref{upper bound for energy} for the details about the choices of $\tilde{\mu}_\hbar$ and $ \tilde{\mathcal{R}}_\hbar$). Using this, we prove that $\mathcal{E}_{\textup{VP}; \min}(M)$ is asymptotically bounded from below by $\mathcal{E}_{\textup{HF}; \min}^\hbar(M)$. 
For the proof, the relativistic Weyl’s law is crucially used to compare the mass and the energy of the auxiliary state and the classical minimizer. 
Conversely, constructing an auxiliary classical state $f_\hbar= \mathbbm{1}_{(\sqrt{1+|p|^2}-1+ \Phi_{\mathcal{Q}_\hbar}(q)<\tilde{\mu}_\hbar')}$,
which is made of the potential function $\Phi_{\mathcal{Q}_\hbar}$ of the quantum minimizer $\mathcal{Q}_\hbar$, and using the  relativistic Weyl’s law, we show the opposite asymptotic inequality.

We remark that compared to the previous work \cite{CHS}, several technical arguments are simplified in this paper, and these simplifications also work in the setting of \cite{CHS}. First of all, we note that for the proof of the minimum energy convergence, we need a version of Weyl’s law with a $\hbar$-parameter family of potentials. 
It is obtained in \cite{CHS} modifying the classical proof in Reed-Simon \cite{RS} by the profile decomposition. 
However, in this article, we instead adopt  the different approach \cite{ELSS, LL}  involving Gaussian coherent states for the relativistic Weyl’s law.
It turns out that this approach is shorter and more direct. Secondly, in \cite{CHS}, rather complicated lemmas are introduced for the desired potential energy convergence. 
However, we found that in the gravitational case, the energy functionals have a coercive structure so that the technical proof in the earlier work can be avoided.

\subsection{Notations}
Throughout this article, the Fourier transform is defined by 
$$\hat{u}(\xi)=\int_{\R^3}u(x)e^{-i x\cdot \xi} dx.$$
For an operator $\gamma: L^2(\mathbb{R}^3)\to L^2(\mathbb{R}^3)$, we denote the operator norm by $\|\gamma\|$ and the Hilbert-Schmidt norm by $\| \gamma \|_{\textup{HS}}:=(\textup{Tr}(\gamma^* \gamma))^\frac12$. 
Given $\phi\in L^2(\mathbb{R}^3)$, $|\phi\rangle\langle\phi|$ denotes the one-particle  projector 
$$|\phi\rangle\langle\phi|:\ u\in L^2(\mathbb{R}^3)\mapsto \int_{\mathbb{R}^3}\phi(x)\overline{\phi(x')}u(x')dx'.$$
 We define 
$$
x_+=\max\{x,0\},\ \  x_-=\min\{x,0\}.
$$

\subsection{Outline of the paper}
The rest of the paper is organized as follows.
In Section 2, we rearrange and reformulate two key inequalities, the Lieb-Thirring inequality and the kinetic interpolation inequality, which fit in our setting. 
The vanishing exchange term estimates  
and the comparison for the quantum and classical critical masses $\mathbf{M}_\text{qm}$ and $\mathbf{M}_\text{cm}$ are also given in Section 2.    
The regularity estimates for the quantum white dwarfs are dealt with in Section 3. 
Section 4 is devoted to provide a version of relativistic Weyl's law with a rate of convergence, which is suitable for our analysis.
We finally prove the main theorem in Section 5 by integrating the aforementioned inequalities, estimates and information.  

\subsection{Acknowledgement}
This work was supported by the New Faculty Startup Fund from Seoul National University. 
This research of the first author was supported by the Basic Science Research Program through the National Research Foundation of Korea (NRF) funded by the Ministry of Science and ICT (RS-2023-00208824 and RS-2023-00219980). This research of the second  author  was supported by the Basic Science Research Program through the National Research Foundation of Korea (NRF)  funded by  the Ministry of Science and ICT  (RS-2023-00213407).
This research of the third author is supported by Basic Science Research Program through the National Research Foundation of Korea (NRF) funded by the Ministry of Science and ICT (NRF-2020R1C1C1A01006415).

\section{Key inequalities and Critical masses}

In this section, we recall preliminary inequalities for the quantum and the kinetic variational problems, and briefly explain why critical masses appear in our variational problems.
 
\subsection{Potential energy estimates}

The kinetic interpolation and the Lieb-Thirring inequalities are employed to control potential energies.

\begin{lemma}[Kinetic interpolation inequality; endpoint case]\label{Kinetic interpolation inequality}
If $0\leq f\leq 1$, then
$$\|\rho_f\|_{L^{4/3}(\mathbb{R}^3)}^{4/3}\lesssim\big\||p|f\big\|_{L^1(\mathbb{R}^3\times\mathbb{R}^3)}.$$
\end{lemma}

\begin{proof}
The density function $\rho_f$ satisfies the following trivial inequality 
$$\rho_f=\int_{|p|\leq R}+\int_{|p|\geq R} f(\cdot,p) dp\lesssim R^3+\frac{1}{R}\big\||p| f(\cdot,p)\big\|_{L_p^1(\mathbb{R}^3)}.$$
Optimizing the right hand side, we obtain $(\rho_f)^{4/3}\lesssim\||p| f(\cdot,p)\|_{L_p^1(\mathbb{R}^3)}$. Then, integrating, we obtain the desired inequality.
\end{proof}

As a quantum analogue of the kinetic interpolation inequality, we have the Lieb-Thirring inequality. Here, for later analysis (see Lemma \ref{regularity lemma}), we give a slightly extended one.

\begin{lemma}[Lieb-Thirring inequality]\label{Lieb-Thirring inequality}
Let $\hbar\in(0,1]$, $s\geq 0$ and $\alpha\in[0,\frac{3}{2})$. If $\gamma$ is a compact self-adjoint operator on $L^2(\mathbb{R}^3)$ and $0\leq |\hbar\nabla|^\alpha\gamma|\hbar\nabla|^\alpha\leq 1$, then 
$$\|\rho_\gamma^\hbar\|_{L^\frac{3+2s-2\alpha}{3-2\alpha}(\mathbb{R}^3)}^{\frac{3+2s-2\alpha}{3-2\alpha}}\lesssim \textup{Tr}^\hbar(|\hbar\nabla|^s\gamma |\hbar\nabla|^s),$$
where the implicit constant is independent of $\hbar$. 
\end{lemma}

\begin{proof}
By scaling $\gamma(x,x')=\frac{1}{\hbar^3}\tilde{\gamma}(\frac{x}{\hbar},\frac{x'}{\hbar})$, we may assume $\hbar=1$. Following the proof in \cite{Sabin}, we apply the Littlewood-Paley inequality in \cite{Sabin} to obtain 
$$\|\rho_\gamma\|_{L^\frac{3+2s-2\alpha}{3-2\alpha}(\mathbb{R}^3)}^{\frac{3+2s-2\alpha}{3-2\alpha}}\lesssim\bigg\|\sum_{N\in 2^{\mathbb{Z}}}\rho_{P_N\gamma P_N}\bigg\|_{L^\frac{3+2s-2\alpha}{3-2\alpha}(\mathbb{R}^3)}^{\frac{3+2s-2\alpha}{3-2\alpha}},$$
where $P_N$ is the Fourier multiplier with symbol $\psi(\frac{\xi}{N})$ such that $\sum_{N\in 2^{\mathbb{Z}}}\psi(\frac{\xi}{N})\equiv1$ on $\mathbb{R}^3$ and $\psi(\xi)$ is supported in $\frac{1}{2}\leq|\xi|\leq4$. We note that as quadratic forms,
$$0\leq P_N\gamma P_N=P_N|\nabla|^{-\alpha}(|\nabla|^{\alpha}\gamma|\nabla|^{\alpha}) |\nabla|^{-\alpha}P_N\lesssim P_N|\nabla|^{-2\alpha}P_N= |\nabla|^{-2\alpha}P_N^2,$$
which implies that
$$0\leq\rho_{P_N\gamma P_N}(x)\leq \rho_{|\nabla|^{-2\alpha}P_N^2}(x)=N^{3-2\alpha}(|\cdot|^{-2\alpha}\psi^2)^\vee(N(x-x'))|_{x=x'}\sim N^{3-2\alpha}.$$
On the other hand, we have $0\leq\rho_{P_N\gamma P_N}\lesssim \frac{1}{N^{2s}}\rho_{|\nabla|^s\gamma |\nabla|^s}$. Thus, it follows that 
$$\sum_{N\in 2^{\mathbb{Z}}}\rho_{P_N\gamma P_N}\lesssim \sum_{N\leq R}N^{3-2\alpha}+\sum_{N\geq R}\frac{1}{N^{2s}}\rho_{|\nabla|^s\gamma |\nabla|^s}\sim R^{3-2\alpha}+\frac{1}{R^{2s}}\rho_{|\nabla|^s\gamma |\nabla|^s}.$$
Thus, optimizing bound as in the proof of Lemma \ref{Kinetic interpolation inequality}, we prove that 
$$\|\rho_\gamma\|_{L^\frac{3+2s-2\alpha}{3-2\alpha}(\mathbb{R}^3)}^{\frac{3+2s-2\alpha}{3-2\alpha}}\lesssim \int_{\mathbb{R}^3}\rho_{|\nabla|^s\gamma |\nabla|^s}dx=\textup{Tr}(|\nabla|^s\gamma |\nabla|^s).$$
\end{proof}

\subsection{Exchange term estimates}

For a compact self-adjoint operator $\gamma$ having a spectral representation $\sum_{j=1}^\infty\lambda_j|\phi_j\rangle\langle\phi_j|$ with an orthonormal set $\{\phi_j\}_{j=1}^\infty$, the exchange term $\textup{X}_{\gamma}^\hbar$ is defined as the integral operator with kernel $\frac{(2\pi\hbar)^3\gamma(x,x')}{|x-x'|}=(2\pi\hbar)^3\sum_{j=1}^\infty\lambda_j\frac{\phi_j(x)\overline{\phi_j(x')}}{|x-x'|}$. By the Fourier transform and the Cauchy-Schwarz inequality, one can see that if $\gamma$ is non-negative, then $0\leq X_{\gamma}^\hbar\leq \frac{1}{|x|}*\rho_{\gamma}^\hbar$ as a quadratic form. 

The following lemma shows that the exchange terms vanishes in the semi-classical limit. 

\begin{lemma}[Exchange term estimates]\label{exchange term est}
If $\gamma$ is a compact self-adjoint operator acting on $L^2(\mathbb{R}^3)$ and $0\leq \gamma\leq 1$, then
$$\begin{aligned}
\|X_\gamma^\hbar\|\leq\|X_\gamma^\hbar\|_{\textup{HS}}&\lesssim \sqrt{\hbar}\Big\{\textup{Tr}^\hbar\big(|\hbar\nabla|\gamma|\hbar\nabla|\big)\Big\}^{\frac{1}{2}},\\
\iint_{\mathbb{R}^3\times\mathbb{R}^3} \frac{|(2\pi\hbar)^3\gamma(x,x')|^2}{|x-x'|}dxdx'&\lesssim\hbar^2\textup{Tr}^\hbar\big(|\hbar\nabla|^{\frac{1}{2}}\gamma|\hbar\nabla|^{\frac{1}{2}}\big).
\end{aligned}$$
\end{lemma}

\begin{proof}
For $s=1,\, 1/2$, by the H\"older and the Sobolev inequalities in the Lorentz norm (see \cite{G} for example), we have
$$\begin{aligned}
\iint \frac{|\gamma(x,x')|^2}{|x-x'|^{2s}}dxdx'&=\bigg\|\frac{\gamma(x,x')}{|x-x'|^{s}}\bigg\|_{L_x^2(\mathbb{R}^3; L_{x'}^2(\mathbb{R}^3))}^2\lesssim\|\gamma(x,x')\|_{L_x^2(\mathbb{R}^3; L_{x'}^{\frac{6}{3-2s},2}(\mathbb{R}^3))}^2\\
&\lesssim \big\||\nabla_{x'}|^s\gamma(x,x')\big\|_{L_x^2(\mathbb{R}^3; L_{x'}^2(\mathbb{R}^3))}^2=\big\|\gamma|\nabla|^s\big\|_{\textup{HS}}^2\\
&\leq\|\sqrt{\gamma}\|^2\big\|\sqrt{\gamma}|\nabla|^s\big\|_{\textup{HS}}^2\leq\textup{Tr}\big(|\nabla|^s\gamma|\nabla|^s\big).
\end{aligned}$$
Then, putting $\hbar$, we obtain the lemma.
\end{proof}

\subsection{Critical masses}
The inequalities in the previous subsection are employed to construct energy minimizers, but this is possible only under the assumption that a given mass is strictly less than the critical mass.

To see this in the kinetic description, we assume that $f\in\mathcal{A}_{\textup{cm}}(M)$. Then, interpolation with the trivial bound $\|\rho_f\|_{L^1(\mathbb{R}^3)}\leq M$ yields 
$$\|\rho_f\|_{L^{6/5}(\mathbb{R}^3)}\lesssim M^{1/3}\big\||p|f\big\|_{L^1(\mathbb{R}^3\times\mathbb{R}^3)}^{1/2},$$
and by the Sobolev inequality, 
\begin{equation}\label{kinetic sharp constant inequality}
\big\||p|f\big\|_{L^1(\mathbb{R}^3\times\mathbb{R}^3)}\geq \frac{K_{\textup{cm}}}{M^{2/3}}\iint_{\mathbb{R}^3\times\mathbb{R}^3}\frac{\rho_f(x)\rho_f(x')}{|x-x'|}dxdx',
\end{equation}
where $K_{\textup{cm}}>0$ is the sharp constant for the above inequality. An important remark is that if $M<\mathbf{M}_{\textup{cm}}=(2K_{\textup{cm}})^{3/2}$, the variational problem $\mathcal{E}_{\textup{VP};\min}(M)$ is well-formulated, since the energy of an admissible distribution $f\in\mathcal{A}_{\textup{cm}}(M)$ is bounded from below;
$$\begin{aligned}
\mathcal{E}_{\textup{VP}}(f)&\geq\big\||p|f\big\|_{L^1(\mathbb{R}^3\times\mathbb{R}^3)}-\|f\|_{L^1(\mathbb{R}^3\times\mathbb{R}^3)}-\frac{M^{2/3}}{2K_{\textup{cm}}}\big\||p|f\big\|_{L^1(\mathbb{R}^3\times\mathbb{R}^3)}\\
&\geq \bigg(1-\frac{M^{2/3}}{2K_{\textup{cm}}}\bigg)\big\||p|f\big\|_{L^1(\mathbb{R}^3\times\mathbb{R}^3)}-M.
\end{aligned}$$

On the other hand, if $\gamma\in\mathcal{A}_{\textup{qm}}(M)$, then by the Sobolev inequality and Lemma \ref{Lieb-Thirring inequality} with $s=\frac{1}{2}$ and $\alpha=0$, there exists the sharp constant $K_{\textup{qm}}>0$ such that 
\begin{equation}\label{quantum sharp constant inequality}
\textup{Tr}^\hbar\big(|\hbar\nabla|^\frac{1}{2}\gamma |\hbar\nabla|^{\frac{1}{2}}|\big)\geq \frac{K_{\textup{qm}}}{M^{2/3}}\iint_{\mathbb{R}^3\times\mathbb{R}^3}\frac{\rho_\gamma^\hbar(x)\rho_\gamma^\hbar(x')}{|x-x'|}dxdx'.
\end{equation}
Thus, it follows that if $M<\mathbf{M}_{\textup{qm}}=(2K_{\textup{qm}})^{3/2}$, then
$$\begin{aligned}
\mathcal{E}_{\textup{HF}}^\hbar(\gamma)\geq \bigg(1-\frac{M^{2/3}}{2K_{\textup{qm}}}\bigg)\textup{Tr}^\hbar\big(|\hbar\nabla|^\frac{1}{2}\gamma |\hbar\nabla|^{\frac{1}{2}}|\big)-M,
\end{aligned}$$
where the exchange energy is dropped since it is non-negative. Note that $\mathbf{M}_{\textup{qm}}$ is independent of $\hbar\in(0,1]$. 

\begin{lemma}[Comparison between the quantum and the classical critical masses]\label{critical mass comparison}
$$\mathbf{M}_{\textup{qm}}=(2K_{\textup{qm}})^{3/2}\leq \mathbf{M}_{\textup{cm}}=(2K_{\textup{cm}})^{3/2}.$$
\end{lemma}

\begin{proof}
For a Schwartz function $f\in\mathcal{S}(\mathbb{R}^3\times\mathbb{R}^3)$, we introduce the operator 
$$\gamma_\hbar=\frac{1}{(2\pi\hbar)^3}\iint |\varphi^\hbar_{(q,p)}\rangle\langle \varphi^\hbar_{(q,p)}| f(q,p)dqdp,$$
where $\varphi^\hbar_{(q,p)}(x)=\frac{1}{(\pi\hbar)^{3/4}}e^{-\frac{|x-q|^2}{2\hbar}}e^{\frac{ip\cdot (x-q)}{\hbar}}$ is a coherent state with $(q,p)\in\mathbb{R}^3\times\mathbb{R}^3$, and insert it into the quantum inequality \eqref{quantum sharp constant inequality}. Then, by direct calculations, taking the limit $\hbar\to0$, one can derive the inequality 
$$\big\||p|f\big\|_{L^1(\mathbb{R}^3\times\mathbb{R}^3)}\geq \frac{K_{\textup{qm}}}{M^{2/3}}\iint_{\mathbb{R}^3\times\mathbb{R}^3}\frac{\rho_f(x)\rho_f(x')}{|x-x'|}dxdx'.$$
Therefore, by density, it follows that $K_{\textup{cm}}\geq K_{\textup{qm}}$.
\end{proof}

\section{Quantum and kinetic white dwarfs and their basic properties}\label{sec: construction of white dwarfs}

The lowest energy states, describing white dwarfs, have been constructed both in the kinetic and the quantum settings (see Theorem \ref{thm: quantum existence} and \ref{thm: classical existence}), employing the inequalities in the previous section. From now on, we assume that $M<\mathbf{M}_{\textup{qm}}$, and let $\mathcal{Q}_\hbar$ be the minimum energy state for the quantum variational problem $\mathcal{E}_{\textup{HF}; \min}^\hbar(M)$, and $\mathcal{Q}$ be that for the kinetic variational problem $\mathcal{E}_{\textup{VP}; \min}(M)$. For notational convenience, we denote
$$\Phi_\hbar:=\Phi_{\mathcal{Q}_\hbar}^\hbar,\quad X_\hbar:=X_{\mathcal{Q}_\hbar}^\hbar\quad\textup{and}\quad \Phi:=\Phi_{\mathcal{Q}}.$$
Then, integrating out the momentum variable for the equation
\begin{equation}\label{kinetic Q self-consistent potential eq}
\mathcal{Q}=\mathbbm{1}_{(\sqrt{1+|p|^2}-1+\Phi(q)\leq\mu)},
\end{equation}
we obtain the potential equation 
\begin{equation}\label{kinetic self-consistent potential eq}
-\Delta\Phi =-\frac{1}{3}\Big\{2(\mu-\Phi)_++(\mu-\Phi)_+^2\Big\}^{\frac{3}{2}}.
\end{equation}
Then, by the elliptic regularity with the kinetic interpolation inequality (Lemma \ref{Kinetic interpolation inequality}), one can show that $\Phi\in L^{3,\infty}(\mathbb{R}^3)\cap C^2(\mathbb{R}^3)$. In addition, a minimizing sequence $\{f^{(n)}\}_{n=1}^\infty$ has a compactness property in the sense that passing to a subsequence and up to translation, $\|\nabla(\frac{1}{|x|}*\rho_{f^{(n)}}-\frac{1}{|x|}*\rho_{\tilde{\mathcal{Q}}})\|_{L^2(\mathbb{R}^3)}\to 0$ for some minimizer $\tilde{\mathcal{Q}}$. For the quantum counterpart, we note that $\mathcal{Q}_\hbar$ obeys the self-consistent equation
\begin{equation}\label{quantum Q self-consistent potential eq}
\mathcal{Q}_\hbar=\mathbbm{1}_{(\sqrt{1-\hbar^2\Delta}-1+\Phi_\hbar+X_\hbar<\mu_\hbar)}+\mathcal{R}_\hbar.
\end{equation}
Moreover, if $\{\gamma^{(n)}\}_{n=1}^\infty$ is a minimizing sequence and up to translation, then passing to a subsequence, $\|\nabla(\frac{1}{|x|}*\rho^\hbar_{\gamma_\hbar}-\frac{1}{|x|}*\rho^\hbar_{\tilde{\mathcal{Q}}_\hbar})\|_{L^2(\mathbb{R}^3)}\to 0$ for some minimizer $\tilde{\mathcal{Q}}_\hbar$.

By construction, a quantum minimizer is a smooth finite-rank operator. To see this, let 
$$\mu_1^\hbar\leq \mu_2^\hbar\leq\mu_3^\hbar\leq\cdots<0$$
denote the negative eigenvalues (counting multiplicities) for $\sqrt{1-\hbar^2\Delta}-1+\Phi_\hbar+X_\hbar$. Then, it is deduced from \eqref{quantum Q self-consistent potential eq} that $\mathcal{Q}_\hbar$ is of the form 
\begin{equation}\label{quantum Q self-consistent potential eq'-1}
\mathcal{Q}_\hbar=\sum_{\mu_j^\hbar<\mu_\hbar}|\phi_j^\hbar\rangle\langle\phi_j^\hbar|+\sum_{\mu_j^\hbar=\mu_\hbar}\lambda_j^\hbar|\phi_j^\hbar\rangle\langle\phi_j^\hbar|\quad \textup{for some }\lambda_j^\hbar\in[0,1],
\end{equation}
where each $\phi_j^\hbar$ solves the eigenvalue equation
\begin{equation}\label{quantum Q self-consistent potential eq'-2}
\big(\sqrt{1-\hbar^2\Delta}-1+\Phi_\hbar+X_\hbar\big)\phi_j^\hbar=\mu_j^\hbar\phi_j^\hbar.
\end{equation}
By the structure \eqref{quantum Q self-consistent potential eq'-1}, $\mathcal{Q}_\hbar$ must have finite rank. Note that the equation \eqref{quantum Q self-consistent potential eq'-2}, together with \eqref{quantum Q self-consistent potential eq'-1}, is simply a finitely many   coupled system of energy sub-critical elliptic PDEs. Standard elliptic regularity theory ensures that all $\phi_j^\hbar$'s are smooth, and so are the total density $\rho_{\mathcal{Q}_\hbar}^\hbar$ and the potential $\Phi_\hbar$. 
We give a simple proof for this in Appendix.   However, obtaining uniform bounds is not so obvious.

In the following lemma, we prove uniform regularity. 

\begin{lemma}[Regularity of quantum white dwarfs]\label{regularity lemma}
For $M\in(0,\mathbf{M}_{\textup{qm}})$, let $\mathcal{Q}_\hbar$ be a minimizer for the variational problem $\mathcal{E}_{\min}^\hbar(M)$ and denote the corresponding potential function by $\Phi_\hbar$. 
Then, there exists some $\alpha \in (0,1)$ such that  
$$\sup_{0<\hbar\leq1}\|\Phi_\hbar\|_{L^{3,\infty}(\mathbb{R}^3)\cap C^{1,\alpha}(\mathbb{R}^3)}<\infty.$$
\end{lemma}

\begin{proof}
By the Sobolev and the Lieb-Thirring (Lemma \ref{Lieb-Thirring inequality}) inequalities, $\|\Phi_\hbar\|_{L^{3,\infty}(\mathbb{R}^3)\cap L^{12}(\mathbb{R}^3)}$ is uniformly bounded. To upgrade regularity, we write from \eqref{quantum Q self-consistent potential eq'-1} that 
$$(2\pi\hbar)^3\sum_{\mu_j^\hbar<\mu_\hbar}(\mu_j^\hbar-\mu_\hbar)^2=\textup{Tr}^\hbar\big((\sqrt{1-\hbar^2\Delta}-1+\Phi_\hbar+X_\hbar-\mu_\hbar)^2\mathcal{Q}_\hbar\big).$$
Note that all quantities in the above identity are finite, because $\mathcal{Q}_\hbar$ is smooth and of finite rank. Moreover, we have
$$
(2\pi\hbar)^3\sum_{\mu_j^\hbar<\mu_\hbar}(\mu_j^\hbar-\mu_\hbar)^2=\mu_\hbar^2(2\pi\hbar)^3\sum_{\mu_j^\hbar<\mu_\hbar}\Big(\frac{\mu_j^\hbar}{\mu_\hbar}-1\Big)^2\le \|\Phi_\hbar\|_{L^\infty(\mathbb{R}^3)}^2\textup{Tr}^\hbar(\mathcal{Q}_\hbar)
$$
and
$0\leq X_\hbar\leq-\Phi_\hbar$ as quadratic forms, and thus   $\|X_\hbar\|\le \|\Phi_\hbar\|= \|\Phi_\hbar\|_{L^\infty(\mathbb{R}^3)}.$
 Hence, by the Cauchy-Schwarz inequality, we obtain
$$\begin{aligned}
\|\Phi_\hbar\|_{L^\infty(\mathbb{R}^3)}^2 M&\geq\textup{Tr}^\hbar\big((1-\hbar^2\Delta)\mathcal{Q}_\hbar\big)-\|-1-\mu_\hbar+\Phi_\hbar+X_\hbar\|^2\textup{Tr}^\hbar(\mathcal{Q}_\hbar)\\
&\quad-2(2\pi\hbar)^3\big\|\sqrt{1-\hbar^2\Delta}\sqrt{\mathcal{Q}_\hbar}\big\|_{\textup{HS}}\|\sqrt{\mathcal{Q}_\hbar}\|_{\textup{HS}}\|-1-\mu_\hbar+\Phi_\hbar+X_\hbar\|\\
&\geq \frac{1}{2}\textup{Tr}^\hbar\big((1-\hbar^2\Delta)\mathcal{Q}_\hbar\big)-9(1+3\|\Phi_\hbar\|_{L^\infty(\mathbb{R}^3)})^2M.
\end{aligned}$$  
However, by interpolation and  the Lieb-Thirring inequality (Lemma \ref{Lieb-Thirring inequality}), we have
$$\|\Phi_\hbar\|_{L^\infty(\mathbb{R}^3)}=\|\tfrac{1}{|x|}*\rho_{\mathcal{Q}_\hbar}^\hbar\|_{L^\infty(\mathbb{R}^3)}\lesssim\|\rho_{\mathcal{Q}_\hbar}^\hbar\|_{L^{\frac{4}{3}}(\mathbb{R}^3)}^{\frac{4}{9}}\|\rho_{\mathcal{Q}_\hbar}^\hbar\|_{L^{\frac{5}{3}}(\mathbb{R}^3)}^{\frac{5}{9}}\lesssim\Big\{\textup{Tr}^\hbar\big((-\hbar^2\Delta)\mathcal{Q}_\hbar\big)\Big\}^{\frac{1}{3}}.$$
Thus, by Young's inequality,  one can show that $\textup{Tr}^\hbar((1-\hbar^2\Delta)\mathcal{Q}_\hbar)$ and $\|\Phi_\hbar\|_{L^\infty(\mathbb{R}^3)}$ are uniformly bounded. Similarly, by the Cauchy-Schwarz inequality and the fact that
\begin{align*}
\|\sqrt{1-\hbar^2\Delta}\mathcal{Q}_\hbar(-1+\Phi_\hbar+X_\hbar)\|^2&\le \|\sqrt{1-\hbar^2\Delta}\sqrt{\mathcal{Q}_\hbar}\sqrt{\mathcal{Q}_\hbar}\|^2(1+2\|\Phi_\hbar\|_{L^\infty(\mathbb{R}^3)})^2\\
&\le \|\sqrt{1-\hbar^2\Delta}\mathcal{Q}_\hbar\sqrt{1-\hbar^2\Delta}\| (1+2\|\Phi_\hbar\|_{L^\infty(\mathbb{R}^3)})^2,
\end{align*}
 one can show that  
$$\begin{aligned}
1\gtrsim(\mu_1^\hbar)^2&=\big\| (\sqrt{1-\hbar^2\Delta}-1+\Phi_\hbar+X_\hbar)\mathcal{Q}_\hbar(\sqrt{1-\hbar^2\Delta}-1+\Phi_\hbar+X_\hbar) \big\|\\
&\geq\frac{1}{2}\|\sqrt{1-\hbar^2\Delta}\mathcal{Q}_\hbar\sqrt{1-\hbar^2\Delta}\|-9(1+2\|\Phi_\hbar\|_{L^\infty(\mathbb{R}^3)})^2M,
\end{aligned}$$  
and so $\|\sqrt{1-\hbar^2\Delta}\mathcal{Q}_\hbar\sqrt{1-\hbar^2\Delta}\|$ is also uniformly bounded. Then, the Lieb-Thirring inequality (Lemma \ref{Lieb-Thirring inequality}) with $\alpha=s=1$ yields that $\|\rho_{\mathcal{Q}_\hbar}^\hbar\|_{L^3(\mathbb{R}^3)}\lesssim 1$.

Repeating the same argument to 
$$(2\pi\hbar)^3\sum_{\mu_j^\hbar<\mu_\hbar}(\mu_j^\hbar-\mu_\hbar)^4=\textup{Tr}^\hbar\big((\sqrt{1-\hbar^2\Delta}-1+\Phi_\hbar+X_\hbar-\mu_\hbar)^4\mathcal{Q}_\hbar\big)$$
and
$$(\mu_j^\hbar-\mu_\hbar)^4=\textup{Tr}^\hbar\big\|(\sqrt{1-\hbar^2\Delta}-1+\Phi_\hbar+X_\hbar-\mu_\hbar)^2\mathcal{Q}_\hbar(\sqrt{1-\hbar^2\Delta}-1+\Phi_\hbar+X_\hbar-\mu_\hbar)^2\big\|,$$
one can show that $\textup{Tr}^\hbar((1-\hbar^2\Delta)^2\mathcal{Q}_\hbar)$ and $\|(1-\hbar^2\Delta)\mathcal{Q}_\hbar(1-\hbar^2\Delta)\|$ are uniformly bounded. Therefore, we obtain that as quadratic forms,
$$0\leq \mathcal{Q}_\hbar\lesssim (1-\hbar^2\Delta)^{-2} $$
so that
\begin{equation}\label{uniform-bdd}
0\leq\rho_{\mathcal{Q}_\hbar}^\hbar\lesssim (2\pi\hbar)^3\rho_{(1-\hbar^2\Delta)^{-2}}\sim 1.
\end{equation} 
Then the elliptic estimate shows that there exists some $ \alpha \in (0,\,1)$ such that $\Phi_\hbar$ is uniformly bounded in $C^{1,\alpha}(\mathbb{R}^3)$.
\end{proof}

\section{Relativistic Weyl's law}

Weyl's law is a fundamental tool for semi-classical analysis as it provides a precise asymptotic formula for eigenvalue counting for a Schr\"odinger operator in terms of the volume of level sets with respect to the corresponding classical local Hamiltonian. Weyl's law has been established in   various settings \cite{ELSS, Frank1, LL},  but we present a version of it is presented in detail (Proposition \ref{Weyl's law}), because we could not find the exact formulation we need.

\subsection{Setup}
We begin this section with some remarks on our setting. Our goal is to develop a Weyl's law, which applies to the quantum minimizer
\begin{equation}\label{self consistent equation for Q}
\mathcal{Q}_\hbar= \mathbbm{1}_{(\sqrt{1-\hbar^2\Delta}-1+\Phi_\hbar+X_\hbar<\mu_\hbar)}+\mathcal{R}_\hbar.
\end{equation}
Hence, we need to handle a family of potentials as well as exchange terms with the parameter which we will take a limit. On the other hand, we may take advantages from that potentials are regular (see Lemma \ref{regularity lemma}) and that only eigenvalues away from the bottom of the continuous spectrum, not all negative eigenvalues, are taken in account (see the formula \eqref{quantum Q self-consistent potential eq}). We also note that the exchange term is not so essential in the semi-classical limit, because it vanishes eventually (see Lemma \ref{exchange term est}).

For this reason, the following assumptions are imposed throughout this section. 
\begin{assumption}[Weyl's law assumptions]\label{Weyl's law assumption}
For a family $\{V_\hbar\}_{\hbar\in(0,1]}$ of potentials and a family $\{A_\hbar\}_{\hbar\in(0,1]}$ of self-adjoint operators on $L^2(\mathbb{R}^3)$, the following hold.
\begin{enumerate}
\item $V_\hbar\not\equiv 0$ is non-positive;
\item $\|V_\hbar\|_{ L^{3,\infty}(\R^3)\cap C^1(\mathbb{R}^3)}$ is bounded uniformly in $\hbar\in(0,1]$;
\item There exists $E_0<0$ such that $\|(V_\hbar-\tfrac{E_0}{2})_-\|_{L^{3/2}(\mathbb{R}^3)}$ is bounded uniformly in $\hbar\in(0,1]$;
\item $\|A_\hbar\|_{\textup{HS}}=O(\sqrt{\hbar})$.
\end{enumerate}
\end{assumption}

\subsection{Eigenvalue counting}
For each $\hbar>0$, we denote the negative eigenvalues of the operator $\sqrt{1-\hbar^2\Delta}-1+V_\hbar+A_\hbar$ in non-decreasing order (counting multiplicities) by
$$E_1^\hbar\le E_2^\hbar\leq E_3^\hbar\leq\cdots<0$$
(see \cite[Theorem XIII.14]{RS} and \cite[Theorem 4.1]{S}). Given the energy level $E<0$, we denote the number of eigenvalues $< E$ by 
$$N^\hbar(E)=N^\hbar(E;V_\hbar, A_\hbar)=\textup{Tr}\big(\mathbbm{1}_{(\sqrt{1-\hbar^2\Delta}-1+V_\hbar+A_\hbar<E)}\big),$$
and define the associated sum
$$\begin{aligned}
S^\hbar(E)=S^\hbar(E;V_\hbar, A_\hbar)&=(2\pi\hbar)^3\sum_{E_j^\hbar< E} (E_j^\hbar-E)\\
&=\textup{Tr}^\hbar\big((\sqrt{1-\hbar^2\Delta}-1+V_\hbar+A_\hbar-E)\mathbbm{1}_{(\sqrt{1-\hbar^2\Delta}-1+V_\hbar+A_\hbar<E)}\big).
\end{aligned}$$
By the Cwikel-Lieb-Rozenblum (CLR) bound below, $(2\pi\hbar)^3N^\hbar(E)$ is uniformly bounded under Assumption \ref{Weyl's law assumption}.
\begin{lemma}[CLR bound]\label{CLR bound}
If $E_0<0$, then
$$(2\pi\hbar)^3N^\hbar(E)\lesssim \frac{1}{(1-(1+\tfrac{E_0}{2})_+^2)^{3/2}}\|(\tfrac{E_0}{2}-V_\hbar)_+\|_{L^3(\mathbb{R}^3)}^3\quad \textup{for all }E\leq E_0.$$
\end{lemma}

\begin{proof}
We employ the relativistic Rumin's inequality \cite{Rumin1, Rumin2}; any self-joint operator $\gamma$ with $0\leq\gamma \leq (\sqrt{1-\hbar^2\Delta}-1-\frac{E}{2})^{-1}$ satisfies 
\begin{equation}\label{relativistic Rumin's inequality}
\textup{Tr}^\hbar\big(\gamma^{1/2}(\sqrt{1-\hbar^2\Delta}-1-\tfrac{E}{2})\gamma^{1/2}\big)\gtrsim \big(1-(1+\tfrac{E}{2})_+^2\big)^{3/4}\|\rho_{\gamma}^\hbar\|_{L^{3/2}(\mathbb{R}^3)}^{3/2}.
\end{equation}
Indeed, the above inequality immediately follows from the generalized Rumin's inequality in Frank \cite[Lemma 2.5]{Frank}, which states in our context that if $0\leq \gamma\leq a(-i\nabla)$ with $a\in L^{3,\infty}(\mathbb{R}^3)$, then
$\textup{Tr} (\gamma^{1/2}a(-i\nabla)\gamma^{1/2})\geq\frac{4\pi^{3/2}}{3\sqrt{3}}\|a\|_{L^{3,\infty}(\mathbb{R}^3)}^{-3/2} \|\rho_{\gamma}\|_{L^{3/2}(\mathbb{R}^3)}^{3/2}$. The inequality \eqref{relativistic Rumin's inequality} is obtained taking 
$a(\xi)=\frac{1}{\sqrt{1+\hbar^2|\xi|^2}-1-\frac{E}{2}}$, because  $\|a\|_{L^{3,\infty}(\mathbb{R}^3)}^3=\frac{4\pi}{3\hbar^3(1-(1+\frac{E}{2})_+^2)^{3/2}}$.

For the CLR bound, we follow the argument of Frank \cite{Frank}. We take $N^\hbar=N^\hbar(E)$ eigenfunctions corresponding to the eigenvalues $\leq E$ for the operator $\sqrt{1-\hbar^2\Delta}-1+V_\hbar+A_\hbar$. Then, by the Gram-Schmidt process on the span of such eigenfunctions, we can construct $\{\psi_j\}_{j=1}^{N^\hbar}$ such that $\langle(\sqrt{1-\hbar^2\Delta}-1-\frac{E}{2}) \psi_j,\psi_k\rangle_{L^2(\mathbb{R}^3)}=\delta_{jk}$, and we obtain the operator $\gamma=\sum_{j=1}^{N^\hbar}|\psi_j\rangle\langle\psi_j|$ satisfying $0\leq\gamma \leq (\sqrt{1-\hbar^2\Delta}-1-\frac{E}{2})^{-1}$ and 
$\textup{Tr}(\gamma^{1/2}(\sqrt{1-\hbar^2\Delta}-1-\frac{E}{2})\gamma^{1/2})=N^\hbar$. Thus, it follows that for $E\leq E_0$,
$$\begin{aligned}
0&\geq \textup{Tr}^\hbar\big(\gamma^{1/2}(\sqrt{1-\hbar^2\Delta}-1+V_\hbar+A_\hbar-E)\gamma^{1/2}\big)\\
&\geq (2\pi\hbar)^3N^\hbar-\int_{\mathbb{R}^3}(\tfrac{E_0}{2}-V_\hbar)_+\rho_\gamma^\hbar dx-\|A_\hbar\|\textup{Tr}^\hbar(\gamma)\\
&\geq (2\pi\hbar)^3 N^\hbar-\|(\tfrac{E_0}{2}-V_\hbar)_+\|_{L^3(\mathbb{R}^3)}\|\rho_\gamma^\hbar\|_{L^{\frac{3}{2}}(\mathbb{R}^3)}-\sqrt{\hbar}(2\pi\hbar)^3 N^\hbar,
\end{aligned}$$
where Lemma \ref{exchange term est} is used for the exchange term. Then, by \eqref{relativistic Rumin's inequality}, the CLR bound follows.
\end{proof}

\subsection{Relativistic Weyl's law}
We prove the main result of this section.
\begin{proposition}[Relativistic Weyl's law]\label{Weyl's law}
Suppose that $\{V_\hbar\}_{\hbar\in(0,1]}$ and $\{A_\hbar\}_{\hbar\in(0,1]}$ satisfy Assumption \ref{Weyl's law assumption} with some $E_0<0$. Then, for $E_\hbar\leq E_0$, we have
\begin{align}
S^\hbar(E_\hbar)&=\iint_{\sqrt{1+|p|^2}-1+V_\hbar(q)<E_\hbar}\big(\sqrt{1+|p|^2}-1+V_\hbar(q)-E_\hbar\big)dqdp+O(\sqrt{\hbar}),\label{Weyl's law; trace}\\
(2\pi\hbar)^3N^\hbar(E_\hbar)&=\big|\big\{(q,p)\in\mathbb{R}^3\times\mathbb{R}^3: \sqrt{1+|p|^2}-1+V_\hbar(q)\leq E_\hbar\big\}\big|+O(\hbar^{1/4}).\label{eq: Weyl's law}
\end{align}
\end{proposition}

For the proof, among several approaches, we follow the argument in \cite{ELSS, LL} involving coherent states
$$\varphi^\hbar_{(q,p)}(x)=\frac{1}{(\pi\hbar)^{3/4}}e^{-\frac{|x-q|^2}{2\hbar}}e^{\frac{ip\cdot (x-q)}{\hbar}},\quad(q,p)\in\mathbb{R}^3\times\mathbb{R}^3.$$

\begin{lemma}[Basic properties of coherent states]\label{basic properties of coherent sates}
$$\begin{aligned}
\int_{\mathbb{R}^3}|\langle \varphi_{(q,p)}^\hbar | u\rangle|^2 dp&=(2\pi\hbar)^3\int \frac{1}{(\pi\hbar)^{3/2}}e^{-\frac{|q-x|^2}{\hbar}}|u(x)|^2dx,\\
\int_{\mathbb{R}^3}|\langle \varphi_{(q,p)}^\hbar| u\rangle|^2 dq&=\hbar^3\int_{\mathbb{R}^3} \frac{1}{(\pi\hbar)^{3/2}}e^{-\frac{|p-\hbar\xi|^2}{\hbar}}|\hat{u}(\xi)|^2d\xi.
\end{aligned}$$
\end{lemma}

\begin{proof}
By direct calculations, we have
$$\begin{aligned}
\int|\langle \varphi_{(q,p)}^\hbar|u\rangle|^2 dp&=\int \iint \frac{1}{(\pi\hbar)^{3/2}}e^{-\frac{|x-q|^2}{2\hbar}}e^{-\frac{|x'-q|^2}{2\hbar}}e^{\frac{ip\cdot (x-x')}{\hbar}} u(x)\overline{u(x')}dxdx'dp\\
&=\iint (4\pi\hbar)^{\frac{3}{2}}e^{-\frac{|x-q|^2}{2\hbar}}e^{-\frac{|x'-q|^2}{2\hbar}}u(x)\overline{u(x')}\delta(x-x')dxdx'\\
&=(2\pi\hbar)^3\int \frac{1}{(\pi\hbar)^{3/2}}e^{-\frac{|q-x|^2}{\hbar}}|u(x)|^2dx.
\end{aligned}$$
On the other hand, by the Plancherel theorem, 
$$\int|\langle \varphi_{(q,p)}^\hbar|u\rangle|^2 dq=\int|(\varphi_{(0,p)}^\hbar* u)(q)|^2 dq=\hbar^3\int \frac{1}{(\pi\hbar)^{3/2}}e^{-\frac{|p-\hbar\xi|^2}{\hbar}}|\hat{u}(\xi)|^2d\xi,$$
where we used $\widehat{\varphi_{(q,p)}^\hbar}(\xi)=(4\pi\hbar)^{3/4}e^{\frac{|p-\hbar\xi|^2}{2\hbar}}e^{-iq\cdot\xi}$.
\end{proof}

\begin{proof}[Proof of Proposition \ref{Weyl's law}]
\noindent \textbf{\underline{Step 1} (Upper bound for $S^\hbar(E_\hbar)$ in \eqref{Weyl's law; trace})} We introduce the operator 
$$\gamma=\frac{1}{(2\pi\hbar)^3}\iint_{\sqrt{1+|p|^2}-1+V_\hbar(q)\leq E_\hbar}|\varphi_{(q,p)}^\hbar\rangle \langle \varphi_{(q,p)}^\hbar|dqdp.$$
Note that $0\le \gamma \le 1$ and $\gamma$ is of trace class, because
$$\begin{aligned}
\textup{Tr}^\hbar(\gamma)&=\big|\big\{(q,p): \sqrt{1+|p|^2}-1+V_\hbar(q)\leq E_\hbar\big\}\big|\\
&=\frac{4\pi}{3}\int \big\{2(E_\hbar-V_\hbar)_++(E_\hbar-V_\hbar)_+^2\big\}^{\frac{3}{2}} dx<\infty
\end{aligned}$$
(see Section \ref{sec: construction of white dwarfs}). Thus, it follows that 
$$\begin{aligned}
S^\hbar(E_\hbar)&\leq(2\pi\hbar)^3\textup{Tr}\big((\sqrt{1-\hbar^2\Delta}-1+V_\hbar+A_\hbar-E_\hbar)\gamma\big)\\
&=\iint_{\sqrt{1+|p|^2}-1+V_\hbar(q)\leq E_\hbar} \langle\varphi_{(q,p)}^\hbar | (\sqrt{1-\hbar^2\Delta}-1+V_\hbar+A_\hbar-E_\hbar)\varphi_{(q,p)}^\hbar\rangle dqdp.
\end{aligned}$$
On the right hand side, direct calculations via the Fourier transform yield
$$\langle\varphi_{(q,p)}^\hbar | \sqrt{1-\hbar^2\Delta}\varphi_{(q,p)}^\hbar\rangle=\int \sqrt{1+|p+\sqrt{\hbar}\xi|^2}\frac{1}{\pi^{3/2}}e^{-|\xi|^2}d\xi=\sqrt{1+|p|^2}+O(\sqrt{\hbar}),$$
while by Assumption \ref{Weyl's law assumption} (2) and (4), 
$$\langle\varphi_{(q,p)}^\hbar | V_\hbar\varphi_{(q,p)}^\hbar\rangle=\int V_\hbar(x)\frac{1}{(\pi\hbar)^{3/2}}e^{-\frac{|x-q|^2}{\hbar}} dx=V_\hbar(q)+O(\sqrt{\hbar})$$
and
$$|\langle\varphi_{(q,p)}^\hbar | A_\hbar\varphi_{(q,p)}^\hbar\rangle|\leq\|A_\hbar\|\|\varphi_{(q,p)}^\hbar\|_{L^2(\mathbb{R}^3)}^2=O(\sqrt{\hbar}).$$
Here, an important remark is that all $O(\sqrt{\hbar})$ terms are independent of $(q,p)$. Therefore, it follows that 
$$\begin{aligned}
&S^\hbar(E_\hbar)-\iint_{\sqrt{1+|p|^2}-1+V_\hbar(q)<E_\hbar} (\sqrt{1+|p|^2}-1+V_\hbar(q)-E_\hbar)dqdp\\
&\lesssim \sqrt{\hbar}\ \big|\big\{(q,p)\in\mathbb{R}^3\times\mathbb{R}^3: \sqrt{1+|p|^2}-1+V_\hbar(q)<E_\hbar\big\}\big|\\
&= \sqrt{\hbar}\cdot\frac{4\pi}{3}\int \big\{2(E_\hbar-V_\hbar)_++(E_\hbar-V_\hbar)_+^2\big\}^{\frac{3}{2}} dx\lesssim  \sqrt{\hbar}.
\end{aligned}$$
\noindent \textbf{\underline{Step 2} (Lower bound for $S^\hbar(E_\hbar)$ in \eqref{Weyl's law; trace})} By the CLR bound (Lemma \ref{CLR bound}), the operator $\sqrt{1-\hbar^2\Delta}-1+V_\hbar+A_\hbar$ has only finitely many (say $N^\hbar$) eigenvalues $<E_\hbar$. Let $\{u_j\}_{j=1}^{N^\hbar}$ be the corresponding orthonormal eigenfunctions. Then, we have 
$$S^\hbar(E_\hbar)=(2\pi\hbar)^3\sum_{j=1}^{N^\hbar}\langle u_j| (\sqrt{1-\hbar^2\Delta}-1+V_\hbar+A_\hbar-E_\hbar)u_j\rangle.$$
Then, by Lemma \ref{CLR bound} and Assumption \ref{Weyl's law assumption} (2) and (4), 
\begin{align*}
S^\hbar(E_\hbar)&=(2\pi\hbar)^3\sum_{j=1}^{N^\hbar}\big\langle u_j| (\sqrt{1-\hbar^2\Delta}-1+R_\hbar[V_\hbar]-E_\hbar)u_j\big\rangle+O(\sqrt{\hbar})\cdot(2\pi\hbar)^3N^\hbar(E_\hbar)\\
&=(2\pi\hbar)^3\sum_{j=1}^{N^\hbar}\big\langle u_j| (\sqrt{1-\hbar^2\Delta}-1+R_\hbar[V_\hbar]-E_\hbar)u_j\big\rangle+O(\sqrt{\hbar}),
\end{align*}
where $R_\hbar[V_\hbar]=\frac{1}{(\pi\hbar)^{3/2}}e^{-\frac{|\cdot|^2}{\hbar}}*V_\hbar$. Note that by Lemma \ref{basic properties of coherent sates}, 
$$\begin{aligned}
&\frac{1}{(2\pi\hbar)^3}\iint \sqrt{1+|p|^2}|\langle \varphi_{(q,p)}^\hbar | u_j\rangle|^2 dqdp\\
&=\frac{1}{(2\pi)^3}\int \left\{\int_{\mathbb{R}^3}(\sqrt{1+|p+\hbar\xi|^2}\frac{1}{(\pi\hbar)^{3/2}}e^{-\frac{|p|^2}{\hbar}}dp\right\}|\hat{u}_j(\xi)|^2d\xi\\
&=\frac{1}{(2\pi)^3}\int \sqrt{1+\hbar^2|\xi|^2}|\hat{u}_j(\xi)|^2d\xi+O(\sqrt{\hbar})=\langle u_j, \sqrt{1-\hbar^2\Delta}u_j\rangle+O(\sqrt{\hbar})
\end{aligned}$$
and
$$\begin{aligned}
\frac{1}{(2\pi\hbar)^3}\iint |\langle \varphi_{(q,p)}^\hbar|u_j\rangle|^2 dqdp&=\|u_j\|_{L^2(\mathbb{R}^3)}^2,\\
\frac{1}{(2\pi\hbar)^3}\iint V_\hbar(q)|\langle \varphi_{(q,p)}^\hbar| u_j\rangle|^2 dqdp&=\int R_\hbar[V_\hbar]|u_j|^2dx.
\end{aligned}$$
Thus, it follows that 
$$S^\hbar(E_\hbar)=\sum_{j=1}^{N^\hbar}\iint \big(\sqrt{1+|p|^2}-1+V_\hbar(q)-E_\hbar\big)|\langle \varphi_{(q,p)}^\hbar| u_j\rangle|^2 dqdp+O(\sqrt{\hbar}).$$
Finally, applying Bessel's inequality $\sum_{j=1}^{N^\hbar}|\langle \varphi_{(q,p)}^\hbar| u_j\rangle|^2\le \|\varphi_{(q,p)}^\hbar\|_{L^2(\mathbb{R}^3)}^2=1$ and the bathtub principle, we conclude that 
$$S^\hbar(E_\hbar)\geq\iint_{\sqrt{1+|p|^2}-1+V_\hbar(q)<E_\hbar} \big(\sqrt{1+|p|^2}-1+V_\hbar(q)-E_\hbar\big)dqdp-O(\sqrt{\hbar}).$$

\noindent \textbf{\underline{Step 3} (Proof of \eqref{eq: Weyl's law})} By the inequality $(E-E_\hbar)_--(E-(E_\hbar+\hbar^{1/4}))_-\geq \hbar^{1/4} {\bf 1}_{\{E<E_\hbar\}}$, where $a_-=\min\{a,0\}$, we have
$$\begin{aligned}
(2\pi\hbar)^3 N^\hbar(E_\hbar)&\leq(2\pi\hbar)^3\sum_{j=1}^\infty \frac{(E_j-E_\hbar)_--(E_j-(E_\hbar+\hbar^{1/4}))_-}{\hbar^{1/4}}=\frac{S^\hbar(E_\hbar)-S^\hbar(E_\hbar+\hbar^{1/4})}{\hbar^{1/4}}.
\end{aligned}$$
Thus, \eqref{Weyl's law; trace} deduces that 
$$\begin{aligned}
(2\pi\hbar)^3 N^\hbar(E_\hbar)&\leq\frac{1}{\hbar^{1/4}}\iint \big(\sqrt{1+|p|^2}-1+V_\hbar(q)-E_\hbar\big)_-\\
&\qquad\qquad-\big(\sqrt{1+|p|^2}-1+V_\hbar(q)-E_\hbar-\hbar^{1/4}\big)_-dqdp+O(\hbar^{1/4})\\
&=\big|\big\{(q,p)\in \mathbb{R}^3\times\mathbb{R}^3: \sqrt{1+|p|^2}-1+V_\hbar(q)<E_\hbar\big\}\big|+O(\hbar^{1/4}).
\end{aligned}$$
The reverse inequality can be proved repeating the above argument but with the inequality $(E-E_\hbar+h^{1/4})_--(E-E_\hbar)_-\leq \hbar^{1/4}{\bf 1}_{\{E<E_\hbar\}}$.
\end{proof}

\section{Proof of the main theorem}

We break the proof into two parts. The first is to show that the quantum minimum energy is asymptotically bounded from above by the classical minimum energy (Lemma \ref{upper bound for energy}), and the second is to prove the reverse inequality (Lemma \ref{lower bound for energy}).

\begin{lemma}[Asymptotic upper bound on the quantum minimum energy]\label{upper bound for energy}
$$\limsup_{\hbar\to 0}\mathcal{E}_{\textup{HF};\min}^\hbar(M)\le\mathcal{E}_{\textup{VP};\min}(M)<0\quad\textup{for all }M\in(0,\mathbf{M}_{\textup{qm}}).$$
\end{lemma}

\begin{proof}
For the proof, we make use of a minimizer $\mathcal{Q}=\mathbbm{1}_{(\sqrt{1+|p|^2}-1+\Phi\leq\mu)}$ for the classical variational problem $\mathcal{E}_{\textup{VP};\min} (M)$ (see Section \ref{sec: construction of white dwarfs}) to construct a quantum state 
$$\gamma_\hbar=\mathbbm{1}_{(\sqrt{1-\hbar^2\Delta}-1+\Phi< \tilde{\mu}_\hbar)}+\tilde{\mathcal{R}}_\hbar,$$
where $0\leq\tilde{\mathcal{R}}_\hbar\leq 1$ is a self-adjoint operator on the eigenspace of $\sqrt{1-\hbar^2\Delta}-1+\Phi$ associated to a negative eigenvalue $\tilde{\mu}_\hbar$. Here, $\tilde{\mu}_\hbar$ and $\tilde{\mathcal{R}}_\hbar$ are chosen so that $\mathcal{M}^\hbar(\gamma_\hbar)=(2\pi\hbar)^3\textup{Tr}(\gamma_\hbar)=M$ and $\gamma_\hbar$ is admissible for the variational problem $\mathcal{E}_{\textup{HF};\min}^\hbar(M)$. Indeed, $\sqrt{1-\hbar^2\Delta}-1+\Phi$ has infinitely many negative eigenvalues, because $\Phi$ is a long-range potential such that if $|x|$ is large, then $\Phi(x)\leq-\frac{c}{|x|}$ for some $c>0$. Hence, such $\gamma_\hbar$, having mass $M$, always exists. We also observe that $\tilde{\mu}_\hbar\to\mu$. Indeed, if $\tilde{\mu}_{\hbar_j}$ converges to a different strictly negative limit point $\tilde{\mu}<0$, Weyl's law (Proposition \ref{Weyl's law}) deduces a contradiction to $\mathcal{M}^\hbar(\gamma_\hbar)=M$. On the other hand, if $\tilde{\mu}_\hbar\to0$, then Weyl's law implies that $M=\mathcal{M}^\hbar(\gamma_\hbar)\geq \textup{Tr}^\hbar(\mathbbm{1}_{(\sqrt{1-\hbar^2\Delta}-1-\frac{c}{|x|}< \epsilon)})=\|\mathbbm{1}_{(\sqrt{1+|p|^2}-1-\frac{c}{|q|}< \epsilon)}\|_{L^1(\mathbb{R}^3\times\mathbb{R}^3)}$ for any $\epsilon>0$. However, since  $\|\mathbbm{1}_{(\sqrt{1+|p|^2}-1-\frac{c}{|q|}< \epsilon)}\|_{L^1(\mathbb{R}^3\times\mathbb{R}^3)}\to \infty$  as $\epsilon\to 0$, taking sufficiently small $\epsilon>0$, we can deduce a contradiction.

We will show that $\mathcal{E}_{\textup{HF}}^\hbar(\gamma_\hbar)\leq\mathcal{E}_{\textup{VP}}(\mathcal{Q})+o_\hbar(1)$,   which immediately implies that $\mathcal{E}_{\textup{HF};\min}^\hbar(M)\le\mathcal{E}_{\textup{VP};\min}(M)+o_\hbar(1)$. First, we claim that the kinetic energy of $\gamma_\hbar$ is uniformly bounded so that by Lemma \ref{exchange term est}, its exchange energy vanishes;
$$\lim_{\hbar\to0}\iint\frac{|(2\pi\hbar)^3\gamma_\hbar(x,x')|^2}{|x-x'|}dxdx'=0.$$
Indeed, by the H\"older inequality and Weyl's law (Proposition \ref{Weyl's law}), the kinetic energy
$$\textup{Tr}^\hbar\big((\sqrt{1-\hbar^2\Delta}-1)\gamma_\hbar\big)=S^\hbar(\tilde{\mu}_\hbar;\Phi, 0)-\int \Phi\rho_{\gamma_\hbar}^\hbar dx+\tilde{\mu}_\hbar \textup{Tr}^\hbar(\gamma_\hbar)$$
is bounded from above by 
$$\begin{aligned}
&S^\hbar(\tilde{\mu}_\hbar;\Phi, 0)+\|\Phi\|_{L^\infty(\mathbb{R}^3)}\int \rho_{\gamma_\hbar}^\hbar dx+\tilde{\mu}_\hbar \textup{Tr}^\hbar(\gamma_\hbar)\\
&=\iint(\sqrt{1+|p|^2}-1+\Phi-\mu)\mathcal{Q}dqdp+\big(\|\Phi\|_{L^\infty(\mathbb{R}^3)}+\mu\big) M+o_\hbar(1)\\
&=\iint(\sqrt{1+|p|^2}-1)\mathcal{Q}dqdp-\iint\frac{\rho_{\mathcal{Q}}(x)\rho_{\mathcal{Q}}(x')}{|x-x'|}dxdx'+\|\Phi\|_{L^\infty(\mathbb{R}^3)} M+o_\hbar(1) 
\end{aligned}$$  (see Section \ref{sec: construction of white dwarfs}).

Next, we reorganize the terms in the energy of $\gamma_\hbar$ in a similar way as 
$$\begin{aligned}
\mathcal{E}_{\textup{HF}}^\hbar(\gamma_\hbar)&=S^\hbar(\tilde{\mu}_\hbar;\Phi, 0)+\tilde{\mu}_\hbar M+\frac{1}{2}\iint\frac{\rho_{\mathcal{Q}}(x)\rho_{\mathcal{Q}}(x')}{|x-x'|}dxdx'\\
&\quad-\frac{1}{2}\iint\frac{(\rho_{\gamma_\hbar}^\hbar-\rho_{\mathcal{Q}})(x)(\rho_{\gamma_\hbar}^\hbar-\rho_{\mathcal{Q}})(x')}{|x-x'|}dxdx'\\
&\quad+\frac{1}{2}\iint \frac{|(2\pi\hbar)^3\gamma_\hbar(x,x')|^2}{|x-x'|}dxdx'.
\end{aligned}$$
For the upper bound on $\mathcal{E}_{\textup{HF}}^\hbar(\gamma_\hbar)$, on the right hand side, the fourth negative term can be dropped and the last exchange energy vanishes by the claim. Hence, by Weyl's law (Proposition \ref{Weyl's law}) again, it follows that  
$$\begin{aligned}
\mathcal{E}_{\textup{HF}}^\hbar(\gamma_\hbar)&\leq\iint(\sqrt{1+|p|^2}-1+\Phi-\mu)\mathcal{Q}dqdp +\mu M\\
&\qquad+\frac{1}{2}\iint\frac{\rho_{\mathcal{Q}}(x)\rho_{\mathcal{Q}}(x')}{|x-x'|}dxdx'+o_\hbar(1)= \mathcal{E}_{\textup{VP}}^\hbar(\mathcal{Q})+o_\hbar(1).
\end{aligned}$$
Finally, to complete the proof, we recall that $\mathcal{E}_{\textup{VP};\min}(M)=\mathcal{E}_{\textup{VP}}(\mathcal{Q})<0$ \cite[Lemma 2.2]{JS}.
\end{proof}

It is shown that the quantum minimum energy has a uniform negative upper bound. From this, some additional useful information for quantum energy minimizers can be obtained. 

\begin{lemma}\label{Q additional properties}
For $M\in(0,\mathbf{M}_{\textup{qm}})$, let $\mathcal{Q}_\hbar$ be a minimizer for the quantum variational problem $\mathcal{E}_{\textup{HF};\min}^\hbar(M)$ satisfying the self-consistent equation $\mathcal{Q}_\hbar=\mathbbm{1}_{(\sqrt{1-\hbar^2\Delta}-1+\Phi_\hbar+X_\hbar<\mu_\hbar)}+\mathcal{R}_\hbar$. Then, the followings hold.
\begin{enumerate}
\item The quantum potential $\Phi_\hbar$ is uniformly of long-range in the sense that for sufficiently small $\delta>0$, there exists $R_\delta>0$, independent of $\hbar$, such that $\Phi_\hbar(x)\leq -\frac{\delta}{2|x|}$ for all $|x|\geq 2R_\delta$, after a suitable translation.
\item The chemical potentials are uniformly away from zero, i.e., $\displaystyle\limsup_{\hbar\to0}\mu_\hbar<0$.
\end{enumerate}
\end{lemma}

\begin{proof}
For (1), it suffices to show that for any small $\delta>0$, there exist $R_\delta\geq1$, independent of $\hbar>0$, such that up to translation, 
$$\liminf_{\hbar\to0}\|\rho_{\mathcal{Q}_\hbar}^\hbar\|_{L^1(|x|\leq R_\delta)}\geq \delta,$$
since it implies that 
$$\Phi_\hbar(x)=-\int \frac{\rho_{\mathcal{Q}_\hbar}^\hbar(x')}{|x-x'|}dx'\leq -\frac{\delta}{2|x|}\quad\textup{when }|x|\geq 2R_\delta.$$
Indeed, for contradiction,   we assume that for any $R>0$, there exists a sequence $\{\hbar_j\}_{j=1}^\infty$ such that $\hbar_j\to0$ and $\sup_{x'\in\mathbb{R}^3}\|\rho_j\|_{L^1(|\cdot-x'|\leq R)}\to0$ where $\rho_j=\rho_{\mathcal{Q}_{\hbar_j}}^{\hbar_j}$. Then, by \eqref{uniform-bdd},
$$\begin{aligned}
\iint \frac{\rho_j(x)\rho_j(x')}{|x-x'|}dxdx'&=\iint_{|x-x'|\le \frac{1}{R}}+\iint_{\frac{1}{R}\leq |x-x'|\le R}+\iint_{|x-x'|\geq R}\frac{\rho_j(x)\rho_j(x')}{|x-x'|}dxdx'\\
&\lesssim \frac{1}{R^2}\|\rho_j\|_{L^\infty(\mathbb{R}^3)}\|\rho_j\|_{L^1(\mathbb{R}^3)}\\
&\quad+R\left\{\sup_{x'\in\mathbb{R}^3}\int_{\frac{1}{R}\leq |x-x'|\le R} \rho_j(x) dx\right\}\|\rho_j\|_{L^1(\mathbb{R}^3)}+\frac{1}{R}\|\rho_j\|_{L^1(\mathbb{R}^3)}^2\\
&\lesssim \frac{1}{R^2}+R\sup_{x'\in\mathbb{R}^3}\int_{\frac{1}{R}\leq |x-x'|\le R} \rho_j(x) dx+\frac{1}{R}.
\end{aligned}$$
Thus, since $R>0$ can be arbitrarily large, it follows that the potential energy vanishes as $j\to\infty$ and $\mathcal{E}_{\textup{HF};\min}^{\hbar_j}(M)=\mathcal{E}_{\hbar_j}(\mathcal{Q}_{\hbar_j})\geq o_j(1)$. However, it contradicts to Lemma \ref{upper bound for energy}.

For (2), we assume that $\mu_{\hbar_j}\to 0$, and take arbitrarily small $\epsilon>0$. Then, for  sufficiently large $j$, we have 
$$M=\textup{Tr}^{\hbar_j}(\mathcal{Q}_{\hbar_j})\geq\textup{Tr}^{\hbar_j}(\mathbbm{1}_{(\sqrt{1-{\hbar_j}^2\Delta}-1+\Phi_{\hbar_j}+X_{\hbar_j}\leq-\epsilon)}).$$
Then, Weyl's law (Proposition \ref{Weyl's law}) and Lemma \ref{Q additional properties} (1) yield
$$\begin{aligned}
M&\geq\big|\big\{(q,p): \sqrt{1+|p|^2}-1+\Phi_{\hbar_j}(q)<-\epsilon\big\}\big|+o_j(1)\\
&=\frac{4\pi}{3}\int \big\{2(-\epsilon-\Phi_{\hbar_j})_++(-\epsilon-\Phi_{\hbar_j})_+^2\big\}^{\frac{3}{2}} dx+o_j(1)\\
&\geq\frac{4\pi}{3}\int \big\{2(\tfrac{\delta}{2|x|}-\epsilon)_++(\tfrac{\delta}{2|x|}-\epsilon)_+^2\big\}^{\frac{3}{2}} dx+o_j(1).
\end{aligned}$$
However, since $\frac{\delta}{2|x|}\notin L^{\frac{3}{2}}(\mathbb{R}^3)$, taking smaller $\epsilon>0$, we can make the lower bound of the above inequality arbitrarily large. It deduces a contradiction.
\end{proof}

Now, we aim to obtain a lower bound for the quantum minimum energy. For this, we switch the role of the quantum and kinetic energy minimization problems in the proof of Lemma \ref{upper bound for energy}. Using a quantum minimizer $\mathcal{Q}_\hbar$ solving the self-consistent equation $\mathcal{Q}_\hbar= \mathbbm{1}_{(\sqrt{1-\hbar^2\Delta}-1+ \Phi_\hbar+X_\hbar<\mu_\hbar)}+\mathcal{R}_\hbar$ (see Section \ref{sec: construction of white dwarfs}), we make an auxiliary kinetic distribution function
\begin{equation}\label{kinetic ansatz}
f_\hbar= \mathbbm{1}_{(\sqrt{1+|p|^2}-1+ \Phi_\hbar(q)<\tilde{\mu}_\hbar)},
\end{equation}
where $\tilde{\mu}_\hbar$ is chosen so that $f_\hbar$ is admissible for the kinetic variational problem $\mathcal{E}_{\textup{VP}; \min}(M)$, that is, $\mathcal{M}(f_\hbar)=M$. Indeed, such $\tilde{\mu}_\hbar<0$ always exists, because by Lemma \ref{Q additional properties} (1), $
\mathcal{M}(\mathbbm{1}_{(\sqrt{1+|p|^2}-1+ \Phi_\hbar(q)<\mu)}) 
$ increases to infinity as $\mu\uparrow 0$.

In addition, $\tilde{\mu}_\hbar$ is sufficiently close to the quantum chemical potential $\mu_\hbar$.

\begin{lemma}\label{lagrange muliplier convergence}
For $\mu_\hbar$ in \eqref{self consistent equation for Q} and $\tilde{\mu}_\hbar$ in \eqref{kinetic ansatz}, we have
$\tilde{\mu}_\hbar-\mu_\hbar\to0$ as $\hbar\to 0$.
\end{lemma}

\begin{proof}
Applying Weyl's law (Proposition \ref{Weyl's law}) to $\textup{Tr}^\hbar(\mathcal{Q}_\hbar)=M$, we obtain
$$
\mathcal{M}(\mathbbm{1}_{(\sqrt{1+|p|^2}-1+ \Phi_\hbar(q)<\mu_\hbar)})=\mathcal{M}(\mathbbm{1}_{(\sqrt{1+|p|^2}-1+ \Phi_\hbar(q)<\tilde{\mu}_\hbar)})+o_\hbar(1).
$$
Consequently, since
\begin{align*}
&\mathcal{M}(\mathbbm{1}_{(\sqrt{1+|p|^2}-1+ \Phi_\hbar(q)<\mu+\epsilon)})-\mathcal{M}(\mathbbm{1}_{(\sqrt{1+|p|^2}-1+ \Phi_\hbar(q)<\mu)})\\
&\ge
 \int_{\{x:\Phi_\hbar\le \mu\}} \big\{2(\mu+\epsilon-\Phi_{\hbar })+(\mu-\Phi_{\hbar })^2\big\}^{\frac{3}{2}}-\big\{2(\mu -\Phi_{\hbar })+(\mu-\Phi_{\hbar })^2\big\}^{\frac{3}{2}} dx\\
&\ge 2\sqrt{2}\epsilon^\frac32|\{x:\Phi_\hbar\le \mu\}|
\end{align*}
for $0<\epsilon<-\mu$, we deduce
 $\tilde{\mu}_\hbar=\mu_\hbar+o_\hbar(1)$. 
\end{proof}

The following is the reverse inequality for the quantum minimum energy. We also prove that the potential function of $f_\hbar$ approximates the quantum potential $\Phi_\hbar$. Indeed, this additional information will be useful to complete the proof of the main theorem.

\begin{lemma}[Asymptotic lower bound on the quantum minimum energy]\label{lower bound for energy}
Let $M\in(0,\mathbf{M}_{\textup{qm}})$. Then,
$$\mathcal{E}_{\textup{VP};\min}(M)\leq\mathcal{E}_{\textup{VP}}(f_\hbar)\leq\liminf_{\hbar\to 0}\bigg(\mathcal{E}_{\textup{HF};\min}^\hbar(M)-\frac{1}{8\pi} \big\|\nabla\big(\Phi_\hbar-(-\tfrac{1}{|x|}*\rho_{f_\hbar})\big)\big\|_{L^2(\mathbb{R}^3)}^2\bigg),$$
where $\Phi_\hbar=\Phi_{\mathcal{Q}_\hbar}^\hbar$ and $\mathcal{Q}_\hbar$ is a minimizer for the quantum variational problem $\mathcal{E}_{\textup{HF};\min}^\hbar(M)$ and $f_\hbar$ is the distribution function given by \eqref{kinetic ansatz}.
\end{lemma}

\begin{proof}
The kinetic distribution $f_\hbar$ has the energy 
$$\begin{aligned}
   \mathcal{E}_{\textup{VP}}(f_\hbar)&=\iint \big(\sqrt{1+|p|^2}-1+\Phi_\hbar(q)- \mu_\hbar\big)f_\hbar dqdp+\mu_\hbar M\\
   &\quad+\frac{1}{2}\iint\frac{\rho_{\mathcal{Q}_\hbar}^\hbar(x)\rho_{\mathcal{Q}_\hbar}^\hbar(x')}{|x-x'|}dxdx' -\frac{1}{2}\iint\frac{(\rho_{f_\hbar}-\rho_{\mathcal{Q}_\hbar}^\hbar)(x)(\rho_{f_\hbar}-\rho_{\mathcal{Q}_\hbar}^\hbar)(x')}{|x-x'|}dxdx'.
\end{aligned}$$
We observe that by Lemma \ref{lagrange muliplier convergence},
\begin{align*}
&\iint \big(\sqrt{1+|p|^2}-1+\Phi_\hbar(q)- \mu_\hbar\big)f_\hbar dqdp\\
&=\iint \big(\sqrt{1+|p|^2}-1+\Phi_\hbar(q)- \mu_\hbar\big)\mathbbm{1}_{(\sqrt{1+|p|^2}-1+ \Phi_\hbar(q)<\mu_\hbar)} dqdp +o_\hbar(1).
\end{align*}
Then, it follows from Weyl's law (Proposition \ref{Weyl's law})  that  
$$\begin{aligned}
\mathcal{E}_{\textup{VP}}(f_\hbar)&=\textup{Tr}^\hbar\big((\sqrt{1-\hbar^2\Delta}-1+\Phi_\hbar+X_\hbar- \mu_\hbar)\mathcal{Q}_\hbar\big)+ \mu_\hbar M\\
&\quad+\frac{1}{2}\iint \frac{\rho_{\mathcal{Q}_\hbar}^\hbar(x)\rho_{\mathcal{Q}_\hbar}^\hbar(x')}{|x-x'|}dxdx'-\frac{1}{8\pi} \big\|\nabla\big(\Phi_\hbar-(-\tfrac{1}{|x|}*\rho_{f_\hbar})\big)\big\|_{L^2(\mathbb{R}^3)}^2+o_\hbar(1).\\
   &= \mathcal{E}_{\textup{HF}}^\hbar(\mathcal{Q}_\hbar)-\frac{1}{8\pi} \big\|\nabla\big(\Phi_\hbar-(-\tfrac{1}{|x|}*\rho_{f_\hbar})\big)\big\|_{L^2(\mathbb{R}^3)}^2+o_\hbar(1).
\end{aligned}$$  
By definition, $\mathcal{E}_{\textup{VP};\min} (M)\leq\mathcal{E}_{\textup{VP}}(f_\hbar)$ and $\mathcal{E}_{\textup{HF}}^\hbar(\mathcal{Q}_\hbar)=\mathcal{E}_{\textup{HF};\min}^\hbar (M)$. Thus, we complete the proof.
\end{proof}

Now, we are ready to prove the main theorem.

\begin{proof}[Proof of Theorem \ref{main-theorem}]
Combining Lemma \ref{upper bound for energy} and Lemma \ref{lower bound for energy}, we obtain that 
$$\begin{aligned}
\mathcal{E}_{\textup{HF};\min}^\hbar(M)&\leq\mathcal{E}_{\textup{VP};\min}(M)+o_\hbar(1)\leq \mathcal{E}_{\textup{VP}}(f_\hbar) +o_\hbar(1)\\
&=\mathcal{E}_{\textup{HF};\min}^\hbar(M)-\frac{1}{8\pi} \|\nabla(\Phi_\hbar-(-\tfrac{1}{|x|}*\rho_{f_\hbar}))\|_{L^2(\mathbb{R}^3)}^2+o_\hbar(1).
\end{aligned}$$
Therefore, it follows that 
\begin{equation}\label{potential convergence to fh}
\|\nabla(\Phi_\hbar-(-\tfrac{1}{|x|}*\rho_{f_\hbar}))\|_{L^2(\mathbb{R}^3)}\to 0
\end{equation}
and for any sequence $\{\hbar_j\}_{j=1}^\infty$ such that $\hbar_j\to0$, $\{f_{\hbar_j}\}_{j=1}^\infty$ is a minimizing sequence for the variational problem $\mathcal{E}_{\textup{VP};\min}(M)$. 
 Hence, by Theorem \ref{thm: classical existence}.$(i)$,  we have the potential convergence $\|\nabla((-\tfrac{1}{|x|}*\rho_{f_{\hbar_j}})-\Phi)\|_{L^2(\mathbb{R}^3)}\to0$ up to a translation. 
Therefore, together with \eqref{potential convergence to fh}, we conclude that $\|\nabla(\Phi_{\hbar_j}-\Phi)\|_{L^2(\mathbb{R}^3)}\to0$. 
 Combining this with Lemma \ref{regularity lemma}, the Sobolev embedding theorem and the uniform estimate \eqref{uniform-bdd}, we have $\| \Phi_{\mathcal{Q}_{\hbar}}  -\Phi_{\mathcal{Q}}\|_{ C(\R^3)} \to 0$ as $\hbar\to 0$. Then, by Lemma \ref{regularity lemma} and the interpolation inequality,
$$
\|\Phi_{\mathcal{Q}_{\hbar}}  -\Phi_{\mathcal{Q}}\|_{C^1(\R^3)}\le \epsilon \|\Phi_{\mathcal{Q}_{\hbar}}  -\Phi_{\mathcal{Q}}\|_{C^{1,\alpha}(\R^3)}+C_\epsilon \|\Phi_{\mathcal{Q}_{\hbar}}  -\Phi_{\mathcal{Q}}\|_{C(\R^3)}\le O(\epsilon)+C_\epsilon \|\Phi_{\mathcal{Q}_{\hbar}}  -\Phi_{\mathcal{Q}}\|_{C(\R^3)},
$$
where $\alpha\in (0,1)$. Thus, sending $\hbar\to 0$ and then $\epsilon\to 0$, we prove Therem \ref{main-theorem}.   
\end{proof}

\appendix

\section{Elliptic regularity for the equations \ref{quantum Q self-consistent potential eq'-2} }
In this appendix, we prove that for fixed $\hbar>0$, the solutions $\phi_j^\hbar$ to the system of elliptic   equations \eqref{quantum Q self-consistent potential eq'-2} are smooth, even though their high Sobolev norms may depends on $\hbar$.
By simple scaling, we may   assume $\hbar = 1$. Then, it is sufficient to prove the following proposition.

\begin{proposition}\label{reg}
Let $K$ be a positive integer and $\{\lambda_j\}_{j=1}^K$, $\{\nu_j\}_{j=1}^K$ be sequences of positive real numbers. 
Let $\{\phi_j\}_{j=1}^K \in H^{1/2}(\R^3)$ be a solution of  the system of equations
\[
\big(\sqrt{1-\Delta}-1+\nu_j\big)\phi_j = \mathcal{N}_j(\phi_1,\dots, \phi_K),  \quad j = 1,\dots, K, 
\]
where
\[
\mathcal{N}_j(\phi_1,\dots, \phi_K) \coloneqq \sum_{k=1}^K  \lambda_k\left(\frac{1}{|x|} * |\phi_k|^2\right)\phi_j -\sum_{k=1}^K\lambda_k\left(\frac{1}{|x|}*(\phi_j\overline{\phi}_k)\right)\phi_k.
\]
Then, for each $j =1,\dots, K$, $\phi_j$ belongs to $H^s(\R^3)$ for every $s \geq 1/2$.
\end{proposition}

We need the following lemma for estimating the nonlinear term $\mathcal{N}_j(\phi_1,\dots, \phi_K)$. We refer to \cite{CHS2} for the proof. 
\begin{lemma}[Trilinear estimate for the Hartree nonlinearity]\label{Hartree nonlinear estimates}
For any $s \geq 1/2$, there exists a positive constant $C_s$ satisfying
\[
\Big\|\Big(\frac{1}{|x|}*(v_1 v_2)\Big)v_3\Big\|_{H^{s}(\mathbb{R}^3)}  \leq C_s\prod_{j=1}^3\|v_j\|_{H^s(\mathbb{R}^3)},\quad \forall v_1, v_2, v_3 \in H^s(\R^3).
\]
\end{lemma}

\begin{proof}[Proof of Proposition \ref{reg}]
We claim that if $\phi_j \in H^s$ for every $j=1,\dots, K$, then $\phi_j \in H^{s+1}$ for every $j=1,\dots, K$.
Since $\{\phi_j\}_{j=1}^K \in H^{1/2}$, we are done if the claim is shown to be true. 

Suppose that $\phi_j \in H^s$ for every $j = 1,\dots, K$ and some $s \geq 1/2$.
We define $\nu \coloneqq \min\{\nu_1, \dots, \nu_K\}$. It is easy to see that there is a small constant $c > 0$ such that
\[
c(|\xi| +1) \leq \sqrt{1+|\xi|^2}-1+\nu.
\]
Then, by combining this and Lemma  \ref{Hartree nonlinear estimates}, we see that for each $j =1, \dots, K$
\[
\begin{aligned}
\|\phi_j\|_{H^{s+1}} &\leq \frac1c\|(\sqrt{1-\Delta}-1+\nu_j)\phi_j\|_{H^s} = \frac1c\|\mathcal{N}_j(\phi_1,\dots, \phi_K)\|_{H^s}  \leq \frac{2C_s}{c}\sum_{k=1}^K\lambda_k\|\phi_k\|_{H^s}^2\|\phi_j\|_{H^s},
\end{aligned}
\]
which proves the claim.

\end{proof}


\begin{thebibliography}{20}

\bibitem{C} S. Chandrasekhar, \emph{The maximum mass of ideal white dwarfs}, Astrophys. J., \textbf{74} (1931), 81–82.

\bibitem{CHS2} W. Choi, Y. Hong and J. Seok, \emph{Optimal convergence rate and regularity of nonrelativistic limit for the nonlinear pseudo-relativistic equations}, J. Funct. Anal. \textbf{274} (2018), no.3, 695--722.

\bibitem{CHS} W. Choi, Y. Hong and J. Seok, \emph{Semi-classical limit of quantum free energy minimizers for the gravitational Hartree equation}, Arch. Ration. Mech. Anal. \textbf{239} (2021), no. 2, 783--829.

\bibitem{Daubechies} I. Daubechies, \emph{An uncertainty principle for fermions with generalized kinetic energy}, Comm. Math. Phys. \textbf{90} (1983), no. 4, 511--520.

\bibitem{ELSS}
 W.D. Evans, R.T. Lewis, H. Siedentop, J.Ph. Solovej,  \emph{ Counting eigenvalues using coherent states with an application to Dirac and Schr\"odinger operators in the semi-classical limit,} Ark. Mat. \textbf{34} (2), 265--283 (1996).


\bibitem{Frank} R. Frank, \emph{Cwikel's theorem and the CLR inequality}, J. Spectr. Theory \textbf{4} (2014), no. 1, 1--21.

\bibitem{Frank1} 
R. Frank, \emph{Eigenvalue bounds for the fractional Laplacian: a review. Recent developments in nonlocal theory,} 210–235, DeGruyter, Berlin, 2018.


\bibitem{FL}
J. Fr\"ohlich and E. Lenzmann,  \emph{Dynamical collapse of white dwarfs in Hartree- and Hartree-Fock theory}, Comm. Math. Phys. \textbf{274} (2007), no. 3, 737--750.


\bibitem{GS}
R. T. Glassey and J. Schaeffer,   \emph{On symmetric solutions of the relativistic Vlasov-Poisson system}, Comm. Math. Phys. \textbf{101}, 459–473 (1985).

\bibitem{G}
L. Grafakos,  \emph{ Classical Fourier Analysis,} Graduate Texts in Mathematics, vol. 249, 2nd edn. Springer, New York (2008)

\bibitem{HR}
H. Hadzic and G. Rein,    \emph{Global existence and nonlinear stability for the relativistic Vlasov-Poisson system in the gravitational case}, Indiana Univ. Math. J. \textbf{56}, 2453--2488 (2007).


\bibitem{HS} C. Hainzl and B. Schlein, \emph{Stellar collapse in the time dependent Hartree-Fock approximation}, Commun. Math. Phys. \textbf{287} (2009), pp. 705--717.

\bibitem{HLLS}
C. Hainzl, E. Lenzmann, M. Lewin and B. Schlein, \emph{On Blowup for Time-Dependent Generalized Hartree-Fock Equations}, Ann. Henri Poincar\'e \textbf{11} 1023--1052(2010).



\bibitem{JS} J. Jang and J. Seok, \emph{Kinetic description of stable white dwarfs}, Kinet. Relat. Models \textbf{15} (2022), no. 4, 605--620.

\bibitem{LMR}
M. Lemou,  F. M\'ehats and P. Rapha\"el,   \emph{Stable ground states for the relativistic gravitational Vlasov-Poisson system}, Comm. Partial Diff. Eq. \textbf{34}(7), 703--721 (2009)


\bibitem{LeLe} E. Lenzmann and M. Lewin, \emph{Minimizers for the Hartree-Fock-Bogoliubov theory of neutron stars and white dwarfs}, Duke Math. J. \textbf{152} (2010), no. 2, 257--315.

\bibitem{LL} E. H. Lieb and M. Loss, \emph{Analysis}, Graduate Studies in Mathematics, Springer  (2001).


\bibitem{LY} E. H. Lieb and H.-T. Yau, \emph{The Chandrasekhar theory of stellar collapse as the limit of quantum mechanics}, Comm. Math. Phys. \textbf{53} (1987),147--174.

\bibitem{LY1}  E.H. Lieb and H.-T. Yau, \emph{A rigorous examination of the Chandrasekhar theory of stellar collapse}, Astrophys. Jour. \textbf{ 323 } (1987), 140--144.

\bibitem{LS}
T. Luo, J. Smoller,  \emph{Nonlinear dynamical stability of Newtonian rotating and non-rotating white dwarfs and rotating supermassive stars,} Comm. Math. Phys. \textbf{284} (2008), no. 2, 425--457


\bibitem{RS} M. Reed and B. Simon, \emph{Methods of Modern Mathematical Physics: Analysis of Operators. Vol. 4}, Academic Press, New York (1978).

\bibitem{Rumin1} M. Rumin, \emph{Spectral density and Sobolev inequalities for pure and mixed states}, Geom. Funct. Anal. \textbf{20} (2010), no. 3, 817--844.

\bibitem{Rumin2} M. Rumin, \emph{Balanced distribution-energy inequalities and related entropy bounds}, Duke Math. J. \textbf{160} (2011), no. 3, 567--597.

\bibitem{Sabin} J. Sabin, \emph{Littlewood-Paley decomposition of operator densities and application to a new proof of the Lieb-Thirring inequality}, Math. Phys. Anal. Geom. \textbf{19} (2016), no. 2, Art. 11, 11 pp. 

\bibitem{S}
B. Simon,  \emph{ Trace ideals and their applications}, London Math. Soc. Lecture Note Ser., \textbf{ 35 }
Cambridge University Press, Cambridge-New York, 1979, viii+134 pp.

\end{thebibliography}
\end{document}